\documentclass[a4paper,12pt,reqno]{amsart}
\usepackage[T1]{fontenc}
\usepackage[utf8]{inputenc}
\usepackage{bera}
\usepackage{amssymb,dsfont,mathrsfs}
\usepackage[margin=1in]{geometry}
\usepackage{enumitem}
\usepackage{pdflscape,afterpage}
\usepackage[bookmarksdepth=2]{hyperref}

\numberwithin{equation}{section}

\newtheorem{theorem}{Theorem}[section]
\newtheorem{lemma}[theorem]{Lemma}

\newtheorem{proposition}[theorem]{Proposition}

\theoremstyle{definition}
\newtheorem{definition}[theorem]{Definition}

\DeclareMathOperator{\re}{Re}
\DeclareMathOperator{\im}{Im}

\DeclareMathOperator{\sign}{sign}

\DeclareMathOperator{\ind}{\mathds{1}}

\newcommand{\op}{L}
\newcommand{\dn}{K}
\newcommand{\sym}{k}
\newcommand{\ii}{I}
\newcommand{\ops}{\mathbb{L}}
\newcommand{\dns}{\mathbb{K}}
\newcommand{\cont}{\mathscr{C}}
\newcommand{\leb}{\mathscr{L}}

\newcommand{\dom}{\mathscr{D}}
\newcommand{\loc}{{\mathrm{loc}}}

\newcommand{\C}{\mathds{C}}
\newcommand{\R}{\mathds{R}}

\newcommand{\ph}{\varphi}

\newcommand{\thet}{\vartheta}

\renewcommand{\le}{\leqslant}
\renewcommand{\ge}{\geqslant}

\usepackage{ifthen}
\newcommand{\formula}[2][nolabel]%
{%
 \ifthenelse{\equal{#1}{nolabel}}%
 {\begin{align*} #2 \end{align*}}%
 {%
  \ifthenelse{\equal{#1}{}}%
  {\begin{align} #2 \end{align}}%
  {\begin{align} \label{#1} \begin{aligned} #2 \end{aligned} \end{align}}%
 }%
}

\newcommand{\ignore}[1]{}

\begin{document}

\title[Harmonic extension technique for non-symmetric operators]{Harmonic extension technique for non-symmetric operators with completely monotone kernels}
\author{Mateusz Kwaśnicki}
\thanks{Work supported by the Polish National Science Centre (NCN) grant no.\@ 2015/19/B/ST1/01457}
\address{Mateusz Kwaśnicki \\ Faculty of Pure and Applied Mathematics \\ Wrocław University of Science and Technology \\ ul. Wybrzeże Wyspiańskiego 27 \\ 50-370 Wrocław, Poland}
\email{mateusz.kwasnicki@pwr.edu.pl}
\date{\today}

\begin{abstract}
We identify a class of non-local integro-differential operators $\dn$ in $\R$ with Dirichlet-to-Neumann maps in the half-plane $\R \times (0, \infty)$ for appropriate elliptic operators $\op$. More precisely, we prove a bijective correspondence between Lévy operators $\dn$ with non-local kernels of the form $\nu(y - x)$, where $\nu(x)$ and $\nu(-x)$ are completely monotone functions on $(0, \infty)$, and elliptic operators $\op = a(y) \partial_{xx} + 2 b(y) \partial_{x y} + \partial_{yy}$. This extends a number of previous results in the area, where symmetric operators have been studied: the classical identification of the Dirichlet-to-Neumann operator for the Laplace operator in $\R \times (0, \infty)$ with $-\sqrt{-\partial_{xx}}$, the square root of one-dimensional Laplace operator; the Caffarelli--Silvestre identification of the Dirichlet-to-Neumann operator for $\nabla \cdot (y^{1 - \alpha} \nabla)$ with $(-\partial_{xx})^{\alpha/2}$ for $\alpha \in (0, 2)$; and the identification of Dirichlet-to-Neumann maps for operators $a(y) \partial_{xx} + \partial_{yy}$ with complete Bernstein functions of $-\partial_{xx}$ due to Mucha and the author. Our results rely on recent extension of Krein's spectral theory of strings by Eckhardt and Kostenko.
\end{abstract}

\maketitle

%
%                            ---------- o ----------
%

\section{Introduction}
\label{sec:intro}

The purpose of this work is to characterise the class of non-local operators $\dn$ that arise as Dirichlet-to-Neumann maps for certain second-order elliptic operators $\op$ in the half-plane $\R \times (0, \infty)$ (or in a strip $\R \times (0, R)$, with the Dirichlet boundary condition at $y = R$). We assume that $\op$ is translation-invariant with respect to the first variable $x$; that is, the coefficients of $\op$ only depend on the second variable $y$. Thus, we consider elliptic equations of the form
\formula[eq:op:g]{
 & a(y) \partial_{xx} u(x, y) + 2 b(y) \partial_{xy} u(x, y) + c(y) \partial_{yy} u(x, y) \\
 & \hspace*{10em} + d(y) \partial_x u(x, y) + e(y) \partial_y u(x, y) = 0 .
}
In general, the Dirichlet-to-Neumann operator associated to equation~\eqref{eq:op:g} is given by
\formula[eq:dn0:def]{
 \dn f(x) & = \lim_{y \to 0^+} \frac{u(x + \tau(y), y) - u(x, 0)}{\sigma^{-1}(y)} \, ,
}
where $u$ is a solution of~\eqref{eq:op:g} with boundary values $u(x, 0) = f(x)$, $\sigma$ is an appropriate \emph{scale function}, and $\tau$ is an appropriate \emph{shearing profile}; we refer to Sections~\ref{sec:red} and~\ref{sec:dis} for further details. Building upon recent extension of Krein's spectral theory of strings due to Eckhardt and Kostenko~\cite{ek16}, we prove that such Dirichlet-to-Neumann maps are of the form
\formula[eq:dn0]{
 \dn f(x) & = \alpha f''(x) + \beta f'(x) - \gamma f(x) \\
 & \qquad + \int_{-\infty}^\infty (f(x + z) - f(x) - f'(x) z \ind_{(-1, 1)}(z)) \nu(z) dz ,
}
where $\alpha \ge 0$, $\beta \in \R$, $\gamma \ge 0$ and $\nu(z)$ and $\nu(-z)$ are completely monotone functions of $z > 0$. Conversely, if we allow for certain irregularities of the coefficients in~\eqref{eq:op:g}, then every operator $\dn$ given by~\eqref{eq:dn0} is the Dirichlet-to-Neumann map for some \emph{general} elliptic equation~\eqref{eq:op:g}. Furthermore, the corresponding equation~\eqref{eq:op:g} is unique, up to some natural transformations.

In fact, we obtain a \emph{bijective} identification of operators of the form~\eqref{eq:dn0} and Dirichlet-to-Neumann operators corresponding to \emph{reduced} elliptic equations of the form
\formula[eq:eq:g]{
 a(dy) \partial_{xx} u(x, y) + 2 b(y) \partial_{xy} u(x, y) + \partial_{yy} u(x, y) & = 0 ,
}
where, for some $R \in (0, \infty]$, $u$ is a sufficiently regular function on $\R \times [0, R)$, $a(dy)$ is a locally finite measure on $[0, R)$, $b(y)$ is a locally square integrable function on $[0, R)$, and $a(dy) - (b(y))^2 dy$ is non-negative. If $R < \infty$, we impose the Dirichlet boundary condition $u(x, R) = 0$ at $y = R$. A rigorous statement of this result is given in Theorem~\ref{thm:main} below, when precise notions of a solution of~\eqref{eq:eq:g} (Definition~\ref{def:harm}) and the corresponding Dirichlet-to-Neumann operator (Definition~\ref{def:dn}) are available. Two variants of Theorem~\ref{thm:main} are provided in Theorem~\ref{thm:main:2} and Proposition~\ref{prop:main:3}, where different classes of reduced equations are considered. A brief explanation how the general equation~\eqref{eq:op:g} can be transformed into the reduced form~\eqref{eq:eq:g} is given in Section~\ref{sec:red}.

The operator $\dn$ given by~\eqref{eq:dn0} is translation invariant, and hence it is a Fourier multiplier: the Fourier transform of $\dn f$ is the product of $\hat{f}$, the Fourier transform of $f$, and $\hat{\dn}$, the symbol of $\dn$; that is, $\widehat{\dn f}(\xi) = \hat{\dn}(\xi) \hat{f}(\xi)$. Operators $\dn$ of the form~\eqref{eq:dn0} are generators of Lévy processes with completely monotone jumps, introduced by L.C.G.~Rogers in~\cite{r83} and revisited recently by the author in~\cite{k19}. The corresponding symbols $\hat{\dn}$ are called \emph{Rogers functions} in~\cite{k19}, and they are closely related to Nevanlinna--Pick functions. For further discussion, see Proposition~\ref{prop:rogers} below, or~\cite{k19}.

Our main result has an appealing probabilistic interpretation: jump processes that arise as boundary traces of two-dimensional diffusions in a half-plane are Lévy processes with completely monotone jumps, and every Lévy process with completely monotone jumps can be realised as a boundary trace in an essentially unique way, up to natural transformations. Here we assume that diffusions are invariant under translations parallel to the boundary of the half-plane. For a detailed discussion, we refer to a companion paper~\cite{k20}.

Before a detailed statement of our results, we briefly discuss existing literature. The classical Dirichlet-to-Neumann operator $\dn$ in the half-plane $\R \times (0, \infty)$ (or, more generally, in half-space $\R^d \times (0, \infty)$) is defined for the Laplace equation $\partial_{xx} u(x, y) + \partial_{yy} u(x, y) = 0$ (for the half-space, we understand that $\partial_{xx}$ denotes the usual $d$-dimensional Laplace operator). This corresponds to $a(y) = 1$ and $b(y) = 0$ in our notation. For over a century it is well-known that $\dn$ is a non-positive definite unbounded operator on $\leb^2(\R)$ (or $\leb^2(\R^d)$), which satisfies $\dn^2 = -\partial_{xx}$. Thus, $\dn = -(-\partial_{xx})^{1/2}$ in the sense of spectral theory.

A similar representation for an arbitrary fractional power $\dn = -(-\partial_{xx})^{\alpha/2}$ (with $\alpha \in (0, 2)$) of the Laplace operator $\partial_{xx}$ is obtained by setting $a(y) = C_\alpha y^{2/\alpha - 2}$ and $b(y) = 0$ for an appropriate constant $C_\alpha$. That is, fractional powers of the Laplace operator are Dirichlet-to-Neumann operators corresponding to the equation
\formula{
 C_\alpha y^{2/\alpha - 2} \partial_{xx} u(x, y) + \partial_{yy} u(x, y) & = 0 .
}
By a simple change of variable $z = c_\alpha y^\alpha$, one can transform this equation into an equation in divergence form
\formula{
 \nabla_{x,z} \cdot (z^{1 - \alpha} \nabla_{x,z} u)(x, z) = 0 ,
}
more suitable for most application and thus more commonly found in literature. This representation of fractional powers of the Laplace operator was studied already in 1960s (see~\cite{mo69,ms65}), and was definitely stated in the above form by Caffarelli and Silvestre in~\cite{cs07}. We refer to Section~10 in the survey article~\cite{g20} for further discussion.

The Caffarelli--Silvestre extension technique has been extended in various directions, which include, among others, replacing $\partial_{xx}$ with a more general operator and studying solutions in more general function spaces (see, for example, \cite{atw18,bmst16,gms13,g20,k17,st10}). A different approach was taken in~\cite{km18}, where $\partial_{yy}$ was replaced by a general elliptic operator $(a(dy))^{-1} \partial_{yy}$ in the half-line. Using Krein's spectral theory of strings, corresponding Dirichlet-to-Neumann operators $\dn$ have been identified with a certain class of Fourier multipliers: the main result of~\cite{km18} asserts that $\dn = -\psi(-\partial_{xx})$ for a \emph{complete Bernstein function} $\psi$, and the correspondence between $a$ and $\psi$ is bijective. With our notation, this corresponds to elliptic equations~\eqref{eq:eq:g} with $b(y) = 0$, and symmetric Dirichlet-to-Neumann operators $\dn$ given by~\eqref{eq:dn0} with $\beta = 0$ and $\nu(-z) = \nu(z)$ in~\eqref{eq:dn}.

Here we would like to mention that the result of~\cite{km18} is a consequence of Krein's spectral theory of strings, developed in the middle of 20th century. In a similar way, as already mentioned above, our main result follows from the recent extension of Krein's theory developed by Eckhardt and Kostenko in~\cite{ek16}.

The relation between the coefficients $a$ and $b$ of the elliptic equation~\eqref{eq:eq:g} and the coefficients $\alpha$, $\beta$, $\gamma$ and $\nu$ of the corresponding Dirichlet-to-Neumann operator $\dn$ given by~\eqref{eq:dn0} is, unfortunately, very inexplicit, and the author is not aware of any results that would link regularity or asymptotics of $a$ and $b$ with similar properties of the parameters of the Dirichlet-to-Neumann operator $\dn$. The picture is not much different in the symmetric case ($b = \beta = 0$ and $\nu(-z) = \nu(z)$) studied in~\cite{km18}, where the rare examples of such results include asymptotic relations between the parameter $a$ and the symbol of $\dn$ in~\cite{cw20,kw10,kw82}; see also~\cite{dm76} for a detailed treatment of Krein's theory, and~\cite{ssv12} for a more recent overview and historical remarks.

Within the probabilistic context, our results are closely related to (singular) integrals of the local time of the Brownian motion. For further discussion, we refer to~\cite{k20}; here we mention the fundamental work of Biane and Yor~\cite{by87}, as well as more recent~\cite{c01}. We also stress that the elliptic operator $\op$ generates a diffusion the half-plane, and the corresponding Dirichlet-to-Neumann operator is the generator of the trace of this diffusion on the boundary. A detailed discussion is again given in~\cite{k20}, and here we only refer to~\cite{bo87,ksv11,k89,m64,m65} for a sample of related research.

Finally, we comment on the more general case, when the coefficients of the elliptic operator $\op$ are allowed to depend on both $x$ and $y$. Although a lot is known about the corresponding Dirichlet-to-Neumann operators when the coefficients are sufficiently regular (see~\cite{gks19} and the references therein for a sample of such results), to the best knowledge of the author, in this generality, a complete description of the class of corresponding Dirichlet-to-Neumann operators is an open problem. This question for \emph{symmetric} elliptic operators of the form $\op u = \nabla_{x,y} \cdot (A(x, y) \nabla_{x,y} u)$ is closely related to the famous Calderón's question whether the conductivity $A(x, y)$ can be reconstructed by measuring the resistance between different parts of the boundary; this is also known as the \emph{electrical impedance tomography} or \emph{electrical resistivity tomography}, and $A(x, y)$ here can be either a scalar or a symmetric matrix. For a solution of Calderón's question in dimension two, we refer to~\cite{n95}; see~\cite{alp16,u13} for a general overview and further references.

No characterisation is known for the class of Dirichlet-to-Neumann operators corresponding to symmetric elliptic equations discussed in the above paragraph. A natural conjecture in the symmetric case is given in~\cite{im98} in terms of a condition on signs of certain determinants, very similar to the concept of total positivity. In the same paper it is proved that this condition is indeed satisfied, given appropriate regularity of the coefficients. We also refer to~\cite{i00} for a result closely related to the extension technique developed in~\cite{km18}.

Noteworthy, a similar question for symmetric elliptic operators on planar graphs with boundaries has been answered by Colin de~Verdière: a complete characterisation of the corresponding discrete Dirichlet-to-Neumann maps is given in~\cite{c94}.

\subsection{Reduction}
\label{sec:red}

Before we rigorously state our main result, Theorem~\ref{thm:main}, we explain how a \emph{general elliptic equation}~\eqref{eq:op:g} can be transformed into a \emph{reduced elliptic equation} of the form~\eqref{eq:eq:g}. In other words, we show that with no loss of generality we may assume that in the general equation~\eqref{eq:op:g} the linear term is missing ($d(y) = e(y) = 0$), and that the coefficient $c(y)$ at $\partial_{yy}$ is equal to $1$. Within this class the correspondence between the associated Dirichlet-to-Neumann operators $\dn$ defined by~\eqref{eq:dn0:def} and operators given by~\eqref{eq:dn0} is bijective.

For simplicity, in this section we ignore completely all regularity issues. These are discussed in detail in the next section for the class of reduced elliptic equations~\eqref{eq:eq:g}, and only briefly in Section~\ref{sec:dis} for other special cases of~\eqref{eq:op:g}.

We begin with a general elliptic equation of the form~\eqref{eq:op:g}. To be specific, we consider the equation $\op_0 u_0 = 0$ for a function $u_0$ defined on $\R \times (0, R_0)$, where the operator $\op_0$ is given by
\formula{
 \op_0 & = a_0(y) \partial_{xx} + 2 b_0(y) \partial_{xy} + c_0(y) \partial_{yy} + d_0(y) \partial_x + e_0(y) \partial_y .
}
The reduction is divided into three steps.

\emph{Step 1.} Our first transformation is \emph{change of scale}. Let $\sigma(y)$ be an increasing solution of the ordinary differential equation $c_0(y) \sigma''(y) + e_0(y) \sigma'(y) = 0$ in $(0, R_0)$, satisfying $\sigma(0^+) = 0$; this function $\sigma$ is unique up to multiplication by a positive constant. Setting $u_0(x, y) = u_1(x, \sigma(y))$, we find that
\formula{
 \op_0 u_0(x, y) & = a_0(y) \partial_{xx} u_1(x, \sigma(y)) + 2 b_0(y) \sigma'(y) \partial_{xy} u_1(x, \sigma(y)) \\
 & \hspace*{3em} + c_0(y) (\sigma'(y))^2 \partial_{yy} u_1(x, \sigma(y)) + d_0(y) \partial_x u_1(x, \sigma(y)) .
}
The right-hand side is equal to $\op_1 u_1(x, \sigma(y))$, where
\formula{
 \op_1 & = a_1(y) \partial_{xx} + 2 b_1(y) \partial_{xy} + c_1(y) \partial_{yy} + d_1(y) \partial_x
}
and
\formula{
 a_1(\sigma(y)) & = a_0(y) , & \qquad b_1(\sigma(y)) & = b_0(y) \sigma'(y) , \\
 c_1(\sigma(y)) & = c_0(y) (\sigma'(y))^2 , & \qquad d_1(\sigma(y)) & = d_0(y) .
}
In particular, $\op_0 u_0(x, y) = 0$ for $(x, y) \in \R \times (0, R_0)$ if and only if $\op_1 u_1(x, y) = 0$ for $(x, y) \in \R \times (0, R_1)$, where $R_1 = \sigma(R_0^-)$.

\smallskip

\emph{Step 2.} In the next stage, we replace $\op_1$ by
\formula{
 \op_2 & = \frac{1}{c_1(y)} \, \op_1 .
}
Clearly, these two operators correspond to the same class of harmonic functions: if we write $u_2(x, y) = u_1(x, y)$ and $R_2 = R_1$, then $\op_1 u_1(x, y) = 0$ for $(x, y) \in \R \times (0, R_1)$ if and only if $\op_2 u_2(x, y) = 0$ for $(x, y) \in \R \times (0, R_2)$. We have
\formula{
 \op_2 & = a_2(y) \partial_{xx} + 2 b_2(y) \partial_{xy} + \partial_{yy} + d_2(y) \partial_x ,
}
where
\formula{
 a_2(y) & = \frac{a_1(y)}{c_1(y)} \, , & \qquad b_2(y) & = \frac{b_1(y)}{c_1(y)} , & \qquad d_2(y) & = \frac{d_1(y)}{c_1(y)} \, .
}

\smallskip

\emph{Step 3.} The final step of reduction is \emph{shearing}. Let $\tau$ be a solution of the ordinary differential equation $\tau''(y) + d_2(y) = 0$ for $y \in (0, R_2)$ satisfying $\tau(0^+) = 0$; this solution is unique up to addition by a linear term. If we set $u_2(x, y) = u_3(x + \tau(y), y)$, then we find that
\formula{
 \op_2 u_2(x, y) & = (a_2(y) + 2 b_2(y) \tau'(y) + (\tau'(y))^2) \partial_{xx} u_3(x + \tau(y), y) \\
 & \hspace*{3em} + 2 (b_2(y) + \tau'(y)) \partial_{xy} u_3(x + \tau(y), y) + \partial_{yy} u_3(x + \tau(y), y) .
}
The right-hand side is now equal to $\op_3 u_3(x + \tau(y), y)$, where
\formula{
 \op_3 & = a_3(y) \partial_{xx} + 2 b_3(y) \partial_{xy} + \partial_{yy} ,
}
and
\formula{
 a_3(y) & = a_2(y) + 2 b_2(y) \tau'(y) + (\tau'(y))^2 , & \qquad b_3(y) & = b_2(y) + \tau'(y) .
}
Once again, $\op_2 u_2(x, y) = 0$ for $(x, y) \in \R \times (0, R_2)$ if and only if $\op_3 u_3(x, y) = 0$ for $(x, y) \in \R \times (0, R_3)$, where $R_3 = R_2$.

\smallskip

By the above identification, we see that $u_0$ is harmonic with respect to $\op_0$ if and only if $u_3$ is harmonic with respect to $\op_3$. Furthermore,
\formula{
 u_0(x, y) & = u_1(x, \sigma(y)) = u_2(x, \sigma(y)) = u_3(x + \tau(\sigma(y)), \sigma(y)) .
}
Recall that we assume that $\sigma(0^+) = 0$ and $\tau(0^+) = 0$. Thus,
\formula{
 \partial_y u_3(x, 0) & = \lim_{y \to 0^+} \frac{u_3(x, y) - u_3(x, 0)}{y} \\
 & = \lim_{y \to 0^+} \frac{u_0(x + \tau(\sigma(y)), \sigma(y)) - u_0(x, 0)}{y} \\
 & = \lim_{z \to 0^+} \frac{u_0(x + \tau(z), z) - u_0(x, 0)}{\sigma^{-1}(z)} \, .
}
Therefore, the rather non-standard definition~\eqref{eq:dn0:def} of the Dirichlet-to-Neumann operator $\dn_0$ associated to the equation $\op_0 u_0 = 0$ agrees with the Dirichlet-to-Neumann operator $\dn_3$ that corresponds to the equation $\op_3 u_3 = 0$ via the usual formula $\dn_3 f(x) = \partial_y u_3(x, 0)$; here $u_0(x, 0) = u_3(x, 0) = f(x)$.

We stress again that in the above reduction we did not discuss the question of regularity of coefficients and solutions, and this is not merely a technicality. We come back briefly to this question in Section~\ref{sec:dis}. 

\subsection{Assumptions and main result}
\label{sec:main}

In the following part we first give a rigorous definition of the class of reduced elliptic equations~\eqref{eq:eq:g} (Definition~\ref{def:op}) and we carefully define the notion of a solution $u$ of~\eqref{eq:eq:g} (Definition~\ref{def:harm}). Then we prove that every boundary value $f$ corresponds to a unique solution $u$ (Proposition~\ref{prop:harm}). This is used to define the Dirichlet-to-Neumann operator (Definition~\ref{def:dn}). Next, we give a rigorous meaning to the non-local operator $\dn$ given by~\eqref{eq:dn0} (Definitions~\ref{def:cm} and~\ref{def:domain}). Only then we are ready to state our main result, Theorem~\ref{thm:main}. Finally, in Proposition~\ref{prop:rogers} we list some equivalent definitions of the class of operators $\dn$ given by~\eqref{eq:dn0}.

We consider the following class $\ops$ of operators $\op$ corresponding to reduced equations of the form~\eqref{eq:eq:g}.

\begin{definition}
\label{def:op}
We say that $\op$ is an operator of class $\ops$ if and only if, formally,
\formula[eq:op]{
 \op & = a(dy) \partial_{xx} + 2 b(y) \partial_{xy} + \partial_{yy}
}
on $\R \times [0, R)$, where
\begin{enumerate}[label={\rm (\alph*)}]
\item $R \in (0, \infty]$;
\item $a(dy)$ is a non-negative, locally finite measure on $[0, R)$;
\item $b(y)$ is a Borel, real-valued function on $[0, R)$ such that $(b(y))^2$ is locally integrable on $[0, R)$;
\item $a(dy) - (b(y))^2 dy$ is non-negative on $[0, R)$.
\end{enumerate}
We say that $\op$ is of class $\ops^\star$ if additionally $a(\{0\}) = 0$.
\end{definition}

We understand formula~\eqref{eq:op} purely formally: it does not define neither the domain, nor the action of $\op$. In Definition~\ref{def:harm} below, a rigorous meaning is given to the equation $\op u = 0$ for $\op \in \ops$.

Note that we do not assume strict ellipticity of $\op$: when $a(dy) = (b(y))^2 dy$ on some interval, then $\op$ becomes degenerate in the corresponding strip.

As usual, in Definition~\ref{def:op} we identify coefficients $b$ which agree almost everywhere. Whenever we say that $\op$ is an operator of class $\ops$, we use $a(dy)$, $b(y)$ and $R$ for the corresponding parameters described in Definition~\ref{def:op}. We additionally denote the auxiliary parameters
\formula{
 \tilde{a}(dy) & = a(dy) - (b(y))^2 dy && \text{and} & B(y) & = \int_0^y b(t) dt
}
for $y \in [0, R)$.

In this general setting the notion of a solution of the equation $\op u = 0$ (or, in other words, a harmonic function for $\op$) requires a careful formulation. Note that the value $\alpha = a(\{0\})$ has no effect on the following definition, and that the definition automatically requires harmonic functions to be sufficiently regular at infinity.

\begin{definition}
\label{def:harm}
For an operator $\op$ of class $\ops$, a Borel function $u(x, y)$ on $\R \times [0, R)$ is said to be harmonic with respect to $\op$ if:
\begin{enumerate}[label=(\alph*)]
\item\label{def:harm:a} for every $y \in [0, R)$ the function $u(\cdot, y)$ is in $\leb^2(\R)$, and it depends continuously (with respect to the $\leb^2(\R)$ norm) on $y \in [0, R)$; if $R = \infty$, then the $\leb^2(\R)$ norm of $u(\cdot, y)$ is assumed to be a bounded function of $y \in [0, \infty)$, while if $R < \infty$, then we additionally require that $u(\cdot, y)$ converges to zero in $\leb^2(\R)$ as $y \to R^-$;
\item\label{def:harm:b} the function $\tilde{u}(x, y) = u(x + B(y), y)$ is weakly differentiable with respect to $y$ on $\R \times (0, R)$, with the weak derivative denoted by $\partial_y \tilde{u}(x, y)$, and $(\partial_y \tilde{u}(x, y))^2$ is integrable over $\R \times (y_1, y_2)$ whenever $0 < y_1 < y_2 < R$;
\item\label{def:harm:c} the equation $\op u(x, y) = 0$ is satisfied in the weak sense in $\R \times (0, R)$.
\end{enumerate}
\end{definition}

The last item of the above definition requires clarification. If $u$ is sufficiently regular, we can use the usual weak (or distributional) formulation of the equation $\op u = 0$, namely, we require that for every smooth, compactly supported function $v(x, y)$ on $\R \times (0, R)$ we have
\formula[eq:harm:def:weak]{
 & \int_{(0, R)} \biggl(\int_{-\infty}^\infty u(x, y) \partial_{xx} v(x, y) dx\biggr) a(dy) \\
 & \qquad - 2 \int_0^R \biggl(\int_{-\infty}^\infty \partial_y u(x, y) \partial_x v(x, y) dx\biggr) b(y) dy \\
 & \qquad\qquad - \int_0^R \biggl(\int_{-\infty}^\infty \partial_y u(x, y) \partial_y v(x, y) dx\biggr) dy = 0 .
}
However, in the general case, $\partial_y u$ may fail to exist: we only know that $\partial_y \tilde{u}$ is well-defined, where $\tilde{u}(x, y) = u(x + B(y), y)$. Therefore, in the general case we understand condition~\ref{def:harm:c} as
\formula[eq:harm:def]{
 & \int_{(0, R)} \biggl(\int_{-\infty}^\infty u(x, y) \partial_{xx} v(x, y) dx\biggr) \tilde{a}(dy) \\
 & \qquad - 2 \int_0^R \biggl(\int_{-\infty}^\infty \partial_y \tilde{u}(x - B(y), y) \partial_x v(x, y) dx\biggr) b(y) dy \\
 & \qquad\qquad + \int_0^R \biggl(\int_{-\infty}^\infty u(x, y) \partial_{yy} v(x, y) dx\biggr) dy = 0
}
for every smooth, compactly supported function $v(x, y)$ on $\R \times (0, R)$. If $u$ is regular enough, it is straightforward to see that conditions~\eqref{eq:harm:def:weak} and~\eqref{eq:harm:def} are equivalent.

We also clarify that the weak differentiability condition~\ref{def:harm:b} for the $\leb^2(\R)$-valued function $y \mapsto \tilde{u}(\cdot, y)$ is understood in the usual way: there is a locally integrable Borel function $\partial_y \tilde{u}(x, y)$ on $\R \times (0, R)$ such that for every smooth, compactly supported function $v$ on $\R \times (0, R)$, we have
\formula[eq:weak]{
 \int_0^R \int_{-\infty}^\infty \partial_y \tilde{u}(x, y) v(x, y) dx dy & = -\int_0^R \int_{-\infty}^\infty \tilde{u}(x, y) \partial_y v(x, y) dx dy .
}

Again we identify all functions $u$ which are equal almost everywhere; however, we always require continuity of the $\leb^2(\R)$-valued function $y \mapsto u(\cdot, y)$.

The following preliminary result is needed for the definition of the Dirichlet-to-Neumann operator.

\begin{proposition}
\label{prop:harm}
Suppose that $\op$ is an operator of class $\ops$. Then for every $f \in \leb^2(\R)$ there is a unique function $u$ harmonic with respect to $\op$ (in the sense of Definition~\ref{def:harm}) such that $u(x, 0) = f(x)$ for almost all $x \in \R$.
\end{proposition}

Proposition~\ref{prop:harm} is proved in Section~\ref{sec:pde}.

\begin{definition}
\label{def:dn}
For an operator $\op$ of class $\ops^\star$, the Dirichlet-to-Neumann operator $\dn$ associated to the equation $\op u = 0$ is an unbounded operator on $\leb^2(\R)$, defined by the formula
\formula[eq:dn:def]{
 \dn f(x) & = \partial_y u(x, 0) = \lim_{y \to 0^+} \frac{u(x, y) - u(x, 0)}{y} \, ,
}
where $u$ is a harmonic function for $\op$ described in Proposition~\ref{prop:harm}, with boundary values $u(x, 0) = f(x)$. Here the limit in~\eqref{eq:dn:def} is understood in the $\leb^2(\R)$ sense, and $f$ is in the domain $\dom(\dn)$ of the operator $\dn$ if and only if $f \in \leb^2(\R)$ and the limit in~\eqref{eq:dn:def} exists.

If $\op$ is an operator of class $\ops$ with $\alpha = a(\{0\}) > 0$, then we use the definition
\formula[eq:dn:def:star]{
 \dn f(x) & = \alpha f''(x) + \partial_y u(x, 0) = \alpha f''(x) + \lim_{y \to 0^+} \frac{u(x, y) - u(x, 0)}{y} \, ,
}
and we say that $f$ is in the domain $\dom(\dn)$ if and only if $f, f', f'' \in \leb^2(\R)$ (with the second derivative understood in the weak sense) and the limit in~\eqref{eq:dn:def:star} exists.
\end{definition}

Our main result identifies Dirichlet-to-Neumann operators associated to elliptic equations $\op u = 0$ for $\op \in \ops$ with the following class of non-local operators.

\begin{definition}
\label{def:cm}
We say that an operator $\dn$ is of class $\dns$ if and only if
\formula[eq:dn]{
 \dn f(x) & = \alpha f''(x) + \beta f'(x) - \gamma f(x) \\
 & \qquad + \int_{-\infty}^\infty (f(x + z) - f(x) - f'(x) z \ind_{(-1, 1)}(z)) \nu(z) dz
}
for every smooth, compactly supported function $f(x)$ on $\R$, where:
\begin{enumerate}[label={\rm (\alph*)}]
\item $\alpha \ge 0$, $\beta \in \R$ and $\gamma \ge 0$;
\item $\nu(z)$ is a real-valued function on $\R \setminus \{0\}$ such that $\nu(z)$ and $\nu(-z)$ are completely monotone functions of $z > 0$, and $\int_{-\infty}^\infty \min\{1, z^2\} \nu(z) dz < \infty$.
\end{enumerate}
We say that $\dn$ is of class $\dns^\star$ if $\alpha = 0$.
\end{definition}

Whenever we consider an operator $\dn$ of class $\dns$, we use the notation $\alpha, \beta, \gamma$ and $\nu(z)$ introduced above. Additionally, we always extend $\dn$ to a closed unbounded operator on $\leb^2(\R)$, as described below.

It is well-known that every operator $\dn$ of class $\dns$ is a Fourier multiplier with symbol
\formula[eq:dn:symbol]{
 \hat{\dn}(\xi) & = -\alpha \xi^2 + i \beta \xi - \gamma + \int_{\R} (e^{i \xi z} - 1 - i \xi z \ind_{(-1, 1)}(z)) \nu(z) dz
}
for $\xi \in \R$; see, for example,~\cite{a04,s99}. By this we mean that if $f$ is a smooth, compactly supported function on $\R$, then the Fourier transform of $\dn f$ is given by $\hat{\dn}(\xi) \hat{f}(\xi)$.

\begin{definition}
\label{def:domain}
Every operator $\dn$ of class $\dns$ is automatically extended to an unbounded operator on $\leb^2(\R)$, with domain
\formula[eq:dn:domain]{
 \dom(\dn) & = \{f \in \leb^2(\R) : \hat{\dn} \cdot \hat{f} \in \leb^2(\R)\} ,
}
and defined by
\formula[eq:dn:fourier]{
 \widehat{\dn f}(\xi) & = \hat{\dn}(\xi) \hat{f}(\xi) .
}
\end{definition}

We are now ready to state our main result.

\begin{theorem}
\label{thm:main}
\begin{enumerate}[label={\textnormal{(\alph*)}}]
\item\label{thm:main:a} If $\op$ is an operator of class $\ops$, then the Dirichlet-to-Neumann operator $\dn$ associated to the equation $\op u = 0$ is an operator of class $\dns$.
\item\label{thm:main:b} Every operator $\dn$ of class $\dns$ is the Dirichlet-to-Neumann operator associated to the equation $\op u = 0$ for a unique operator $\op$ of class $\ops$.
\end{enumerate}
\end{theorem}

Theorem~\ref{thm:main} is proved in Section~\ref{sec:pde}. Here we observe that it is sufficient to prove Theorem~\ref{thm:main} for classes $\ops^\star$ and $\dns^\star$ rather than $\ops$ and $\dns$. Indeed, suppose that $\op$ is an operator of class $\ops$ such that $\alpha = a(\{0\})$, and let $\op^\star$ be the corresponding operator of class $\ops^\star$, obtained by replacing $a(dy)$ by $a^\star(dy) = \ind_{(0, \infty)}(y) a(dy)$. The operators $\op$ and $\op^\star$ share the same class of harmonic functions. Thus, if $\dn$ and $\dn^\star$ are the corresponding Dirichlet-to-Neumann operators, then $\dn f = \dn^\star f + \alpha f''$ (see Definition~\ref{def:dn}).

We note that very few explicit pairs of associated operators $\op$ and $\dn$ are known; see Section~\ref{sec:ex} for examples and further discussion. We also remark that if $f$ is in the domain of $\dn$ and $u$ is the harmonic extension of $f$, then the weak derivative $\partial_y u$ is well-defined, and formula~\eqref{eq:dn:def} in the definition of the Dirichlet-to-Neumann operator (Definition~\ref{def:dn}) can be equivalently written as
\formula[eq:dn:def:alt]{
 \dn f(x) & = \lim_{y \to 0^+} \partial_y u(x, y) ,
}
with the limit in $\leb^2(\R)$. This follows from Theorem~\ref{thm:ode}\ref{thm:ode:c} and Lemma~\ref{lem:ext:1} by an argument used in the proof of Theorem~4.3 in~\cite{km18}; we omit the details.

In~\cite{k19}, Fourier symbols $\hat{\dn}$ of operators of class $\dns$ are called \emph{Rogers functions}, and a number of equivalent characterisations of this class of functions is given therein. For completeness, we list them in the following statement.

\begin{proposition}[Theorem~3.3 in~\cite{k19}]
\label{prop:rogers}
Suppose that $\sym(\xi)$ is a continuous function on $\R$, satisfying $\sym(-\xi) = \overline{\sym(\xi)}$ for all $\xi \in \R$. The following conditions are equivalent:
\begin{enumerate}[label=\rm (\alph*)]
\item
$-\sym$ is the Fourier symbol of some operator $\dn$ of class $\dns$, that is, $\sym(\xi) = -\hat{\dn}(\xi)$ for all $\xi \in \R$, with $\hat{\dn}$ given by~\eqref{eq:dn:symbol};
\item
for all $\xi \in \R$ we have
\formula[eq:dn:symbol:int]{
 \sym(\xi) & = \alpha \xi^2 - i \check{\beta} \xi + \gamma + \frac{1}{\pi} \int_{\R \setminus \{0\}} \biggl(\frac{\xi}{\xi + i s} + \frac{i \xi \sign s}{1 + |s|}\biggr) \frac{\mu(ds)}{|s|}
}
for some $\check{\beta} \in \R$ and some non-negative measure $\mu$ on $\R \setminus \{0\}$ such that $\int_{\R \setminus \{0\}} \min\{|s|^{-1}, |s|^{-3}\} \mu(ds) < \infty$;
\item\label{it:r:d}
either $\sym(\xi) = 0$ for all $\xi \in \R$ or for all $\xi \in \R$ we have
\formula[eq:dn:symbol:exp]{
 \sym(\xi) & = c \exp\biggl(\frac{1}{\pi} \int_{-\infty}^\infty \biggl(\frac{\xi}{\xi + i s} - \frac{1}{1 + |s|}\biggr) \frac{\thet(s)}{|s|} \, ds\biggr)
}
for some $c > 0$ and some Borel function $\thet$ on $\R$ with values in $[0, \pi]$;
\item
$\sym$ extends to a holomorphic function in the right complex half-plane $\{\xi \in \C : \re \xi > 0\}$ and $\re (\sym(\xi) / \xi) \ge 0$ whenever $\re \xi > 0$ (that is, $\sym(\xi) / \xi$ is a \emph{Nevanlinna--Pick function} in the right complex half-plane).
\end{enumerate}
\end{proposition}

\subsection{Variants}
\label{sec:dis}

Our main result is stated for the reduced elliptic equation $\op u = 0$, with operator $\op$ of the form
\formula[eq:op:1]{
 \op & = a(dy) \partial_{xx} + 2 b(y) \partial_{xy} + \partial_{yy}
}
in $\R \times (0, R)$, where $R \in (0, \infty]$, $a(dy)$ is a non-negative, locally finite measure, and $b(y)$ is a real-valued function such that $(b(y))^2$ is locally integrable and $a(dy) - (b(y))^2 dy \ge 0$. We choose this variant, because it leads to relatively few technical difficulties, and it is well-suited for a probabilistic interpretation. However, various reformulations of our result are possible, two of which are discussed below. More precisely, first we rephrase our main result for the operators of the form
\formula[eq:op:2]{
 \tilde{\op} & = \tilde{a}(dy) \partial_{xx} + \partial_{yy} + \tilde{d}(y) \partial_x ,
}
and then we specialise our theorem to the class of operators
\formula[eq:op:3]{
 \dot{\op} & = \partial_{xx} + (\dot{a}(y))^{-2} (\dot{a}(y) \partial_y - \dot{b}(y) \partial_x) (\dot{a}(y) \partial_y + \dot{b}(y) \partial_x) ,
}
for appropriate coefficients $\tilde{a}$, $\tilde{d}$, $\dot{a}$ and $\dot{b}$. We will refer to the operator $\op$ of the form~\eqref{eq:op:1}, or the corresponding reduced elliptic equation $\op u = 0$, as an operator or an equation in the \emph{standard form}. Similarly, the terms \emph{Eckhardt--Kostenko form}, and \emph{divergence-like form}, will be used in reference to operators $\tilde{\op}$ of the form~\eqref{eq:op:2}, and operators $\dot{\op}$ of the form~\eqref{eq:op:3}, respectively.

Let us stress that, in principle, it is possible to reverse completely the reduction in Section~\ref{sec:red} and state a result for general equations of the form~\eqref{eq:op:g}. However, a complete description of the class of coefficients $a, b, c, d, e$, for which the corresponding Dirichlet-to-Neumann operator is well-defined, is somewhat problematic. Additionally, one loses the bijective correspondence between coefficients and Dirichlet-to-Neumann operators. For these reasons, we take a different perspective, and we focus on operators given by~\eqref{eq:op:2} and~\eqref{eq:op:3}.

\subsubsection{Eckhardt--Kostenko form}

With $\op$ and $\tilde{\op}$ defined by~\eqref{eq:op:1} and~\eqref{eq:op:2}, the equations $\op u = 0$ and $\tilde{\op} \tilde{u} = 0$ are found to be equivalent by choosing another shearing in Step~3 of reduction. Indeed, let us define
\formula[eq:var12]{
 \tilde{a}(dy) & = a(dy) - (b(y))^2 dy , & \qquad \tilde{d}(y) & = -b'(y) ,
}
and let $u$ and $\tilde{u}$ be related one to the other by the formula
\formula{
 \tilde{u}(x, y) & = u(x + B(y), y) , & \qquad \text{with} &\qquad & B(y) & = \int_0^y b(s) ds .
}
Given enough regularity of $a$, $b$ and $u$, it is now straightforward to show that $\op u = 0$ in $\R \times (0, R)$ if and only if $\tilde{\op} \tilde{u} = 0$ in $\R \times (0, R)$. In the general case, however, some care is needed, as it was the case with Definition~\ref{def:harm}: $\tilde{d}$ is the derivative of an arbitrary locally square-integrable function.

\begin{definition}\label{def:harm:2}
Suppose that $R \in (0, \infty)$, $\tilde{a}$ is a locally finite, non-negative measure on $[0, R)$, and $\tilde{d}$ is the distributional derivative of a locally square-integrable function on $[0, R)$. We say that a function $\tilde{u}(x, y)$ is harmonic with respect to the operator $\tilde{\op}$ given by~\eqref{eq:op:2} if:
\begin{enumerate}[label=(\alph*)]
\item\label{def:harm:2:a} for every $y \in [0, R)$ the function $\tilde{u}(\cdot, y)$ is in $\leb^2(\R)$, and it depends continuously (with respect to the $\leb^2(\R)$ norm) on $y \in [0, R)$; if $R = \infty$, then the $\leb^2(\R)$ norm of $\tilde{u}(\cdot, y)$ is assumed to be a bounded function of $y \in [0, \infty)$, while if $R < \infty$, then we additionally require that $\tilde{u}(\cdot, y)$ converges to zero in $\leb^2(\R)$ as $y \to R^-$;
\item\label{def:harm:2:b} the function $\tilde{u}(x, y)$ is weakly differentiable on $\R \times (0, R)$ with respect to $y$, and $(\partial_y \tilde{u}(\cdot, y))^2$ is integrable over $\R \times (y_1, y_2)$ whenever $0 < y_1 < y_2 < R$;
\item\label{def:harm:2:c} the equation $\tilde{\op} \tilde{u} = 0$ is satisfied in the weak sense in $\R \times (0, R)$, that is, for every smooth, compactly supported function $v$ on $\R \times (0, R)$,
\formula[eq:harm:2:def]{
 & \int_{(0, R)} \biggl(\int_{-\infty}^\infty \tilde{u}(x, y) \partial_{xx} v(x, y) dx\biggr) \tilde{a}(dy) \\
 & \qquad + \int_0^R \biggl(\int_{-\infty}^\infty \tilde{u}(x, y) \partial_{yy} v(x, y) dx\biggr) dy \\
 & \qquad\qquad - \int_0^R \biggl(\int_{-\infty}^\infty \tilde{u}(x, y) \partial_x v(x, y) dx\biggr) \tilde{d}(y) dy = 0 .
}
\end{enumerate}
\end{definition}

We clarify that if $\tilde{d} = -b'$, then the last integral in~\eqref{eq:harm:2:def} should be understood as
\formula{
 \int_0^R \biggl(\int_{-\infty}^\infty \bigl(\partial_y \tilde{u}(x, y) \partial_x v(x, y) + \tilde{u}(x, y) \partial_{xy} v(x, y) \bigr) dx\biggr) b(y) dy ,
}
and in particular this is why weak differentiability of $\tilde{u}$ with respect to $y$ is needed.

It is somewhat technical, but relatively straightforward to prove that $u$ is harmonic with respect to $\op$ in the sense of Definition~\ref{def:harm} if and only if $\tilde{u}$ is harmonic with respect to $\tilde{\op}$ in the sense of Definition~\ref{def:harm:2}. In fact, the only difficulty lies in the proof that conditions~\eqref{eq:harm:def} and~\eqref{eq:harm:2:def} are equivalent. We omit the details.

If $\tilde{u}(x, y) = u(x + B(y), y)$ is a harmonic function for $\tilde{\op}$ with boundary values $f(x) = \tilde{u}(x, 0)$, then, according to Definition~\ref{def:dn}, the corresponding Dirichlet-to-Neumann operator $\dn$ is given by
\formula[eq:dn:2]{
 \dn f(x) & = \alpha f''(x) + \lim_{y \to 0^+} \frac{\tilde{u}(x - B(y), y) - \tilde{u}(x, 0)}{y} \, ,
}
with $\alpha = \tilde{a}(\{0\})$ and with the limit in $\leb^2(\R)$. Note that given only $\tilde{\op}$ (that is, the coefficients $\tilde{a}$ and $\tilde{d}$), there is some ambiguity in the above definition: the function $B$ is defined up to addition by a linear term only, and thus $\dn f$ is only defined up to addition by a first-order term $C f'$ for some $C \in \R$.

As an immediate corollary of Theorem~\ref{thm:main}, we obtain the following result.

\begin{theorem}
\label{thm:main:2}
\begin{enumerate}[label={\textnormal{(\alph*)}}]
\item Under the assumptions listed in Definition~\ref{def:harm:2}, the Dirichlet-to-Neumann operator $\dn$ associated to the equation $\tilde{\op} \tilde{u} = 0$, with $\tilde{\op}$ given by~\eqref{eq:op:2}, is an operator of class $\dns$.
\item Every operator $\dn$ of class $\dns$ is the Dirichlet-to-Neumann operator associated to the equation $\tilde{\op} \tilde{u} = 0$ for a unique triplet of parameters $R$, $\tilde{a}$ and $\tilde{d} = -b'$ satisfying the conditions listed in Definition~\ref{def:harm:2}.
\end{enumerate}
\end{theorem}

Compared to the equation $\op u = 0$ in standard form, studied in Section~\ref{sec:main}, the Eckhardt--Kostenko form $\tilde{\op} \tilde{u} = 0$ is much more closely related to the ODE studied in~\cite{ek16}; see Section~\ref{sec:ode} for further discussion. Additionally, the definition of a solution of the elliptic equation $\tilde{\op} \tilde{u} = 0$ is somewhat simpler. On the other hand, the Eckhardt--Kostenko form presents a number of additional technical difficulties. First of all, one has to work with distributional derivatives of square-integrable functions, that is, with elements of the Sobolev space $H^{-1}_\loc([0, R))$ of negative index; again see Section~\ref{sec:ode}. Furthermore, the definition~\eqref{eq:dn:2} of the Dirichlet-to-Neumann operator is less natural for the Eckhardt--Kostenko form. In fact, as described above, formula~\eqref{eq:dn:2} is ambiguous: it depends on the function $B$, which is not uniquely determined by the coefficient $\tilde{d}$ (this is the reason why we write $\tilde{d} = -b'$ rather than simply $\tilde{d}$ in Theorem~\ref{thm:main:2}). Finally, the equation in standard form turns out to be more convenient than the Eckhardt--Kostenko variant in probabilistic applications, to be discussed in~\cite{k20}.

With the above arguments in mind, in this article we focus on the standard form~\eqref{eq:op:1} considered in Section~\ref{sec:main}, and we limit our discussion of the equation~\eqref{eq:op:2} in Eckhardt--Kostenko form to this section. Note, however, that finding an operator $\op$ which corresponds to a given operator $\tilde{\op}$, or vice versa, presents no difficulties; see formula~\eqref{eq:var12}.

\subsubsection{Divergence-like form}

We now move to the variant given by~\eqref{eq:op:3}. In the symmetric case (corresponding to $b(y) = 0$), discussed in detail in~\cite{km18}, it is often convenient to work with the equation in the divergence form: $\nabla_{x,y} \cdot (\dot{a}(y) \nabla_{x,y} \dot{u}) = 0$, rather than the standard form: $a(dy) \partial_{xx} u + \partial_{yy} u = 0$. Both equations are equivalent by an appropriate change of scale, which corresponds to a different choice of $\sigma$ in Step~1 of the reduction. The equation in the divergence form, however, is less general: not every measure $a(dy)$ corresponds to some coefficient $\dot{a}(y)$. We refer to~\cite{km18} for a detailed discussion.

Below we implement a similar strategy in the non-symmetric case, and again we need to impose additional restrictions on the coefficients $a(dy)$ and $b(y)$; in other words, this approach leads to the representation as Dirichlet-to-Neumann operators for a class of operators $\dn$ strictly smaller than $\dns$.

We study the elliptic equation $\dot{\op} \dot{u} = 0$, where $\dot{\op}$ is given by~\eqref{eq:op:3}. More precisely, we consider the equation
\formula[eq:eq:3]{
 (\dot{a}(y))^2 \partial_{xx} \dot{u}(x, y) + (\dot{a}(y) \partial_y - \dot{b}(y) \partial_x) (\dot{a}(y) \partial_y + \dot{b}(y) \partial_x) \dot{u}(x, y) & = 0 ,
}
which, strictly speaking, corresponds to the equation $(\dot{a}(y))^2 \dot{\op} \dot{u}(x, y) = 0$ with the notation of~\eqref{eq:op:3}. Given enough regularity of the coefficients, equation~\eqref{eq:eq:3} takes form
\formula{
 & \bigl((\dot{a}(y))^2 - (\dot{b}(y))^2\bigr) \partial_{xx} \dot{u}(x, y) + (\dot{a}(y))^2 \partial_{yy} \dot{u}(x, y) \\
 & \hspace*{7em} + \dot{a}(y) \dot{b}'(y) \partial_x u(x, y) + \dot{a}(y) \dot{a}'(y) \partial_y u(x, y) = 0 .
}
This equation again corresponds to an equation $\op u = 0$ for an appropriate operator $\op$ in the standard form~\eqref{eq:op:1}. Before we discuss this relation in detail, however, let us first give a rigorous meaning to~\eqref{eq:eq:3}.

\begin{definition}\label{def:harm:3}
Let $\dot{R} \in (0, \infty)$, and suppose that $\dot{a}$ and $\dot{b}$ are functions on $[0, \dot{R})$ such that $\dot{a}$ and $1 / \dot{a}$ are locally integrable on $[0, \dot{R})$, and $|\dot{b}(y)| \le \dot{a}(y)$ for all $y \in [0, \dot{R})$. We say that $\dot{u}$ is a harmonic function for the operator $\dot{\op}$ given by~\eqref{eq:op:3} if:
\begin{enumerate}[label=(\alph*)]
\item\label{def:harm:3:a} for every $y \in [0, \dot{R})$ the function $\dot{u}(\cdot, y)$ is in $\leb^2(\R)$, and it depends continuously (with respect to the $\leb^2(\R)$ norm) on $y \in [0, \dot{R})$; if $1 / \dot{a}$ is not integrable over $[0, \dot{R})$, then the $\leb^2(\R)$ norm of $\dot{u}(\cdot, y)$ is assumed to be a bounded function of $y \in [0, \dot{R})$, while if the integral of $1 / \dot{a}$ is finite, then we additionally require that $\dot{u}(\cdot, y)$ converges to zero in $\leb^2(\R)$ as $y \to \dot{R}^-$;
\item\label{def:harm:3:b} the function $\dot{u}(x, y)$ is weakly differentiable on $\R \times (0, \dot{R})$ with respect to $y$, and $\dot{a}(y) (\partial_y \dot{u}(x, y))^2$ is integrable over $\R \times (y_1, y_2)$ whenever $0 < y_1 < y_2 < \dot{R}$;
\item\label{def:harm:3:c} the equation $\dot{\op} \dot{u} = 0$ is satisfied in the weak sense in $\R \times (0, \dot{R})$, with weight $(a(y))^{-1}$; that is, for every suitable test function $\dot{v}$ on $\R \times (0, \dot{R})$, we have
\formula[eq:harm:3:def]{
 & \int_0^{\dot{R}} \biggl( \int_{-\infty}^\infty \dot{u}(x, y) \partial_{xx} \dot{v}(x, y) dx \biggr) \dot{a}(y) dy \\
 & \qquad + \int_0^{\dot{R}} \biggl( \int_{-\infty}^\infty \dot{u}(x, y) \bigl(\dot{a}(y) \partial_{x y} \dot{v}(x, y) - \dot{b}(y) \partial_{xx} \dot{v}(x, y)\bigr) dx \biggr) \frac{\dot{b}(y)}{\dot{a}(y)} \, dy \\
 & \qquad\qquad - \int_0^{\dot{R}} \biggl( \int_{-\infty}^\infty \partial_y \dot{u}(x, y) \bigl(\dot{a}(y) \partial_y \dot{v}(x, y) - \dot{b}(y) \partial_x \dot{v}(x, y)\bigr) dx \biggr) dy = 0 .
}
\end{enumerate}
\end{definition}

By a \emph{suitable test function} in item~\ref{def:harm:3:c} we understand a compactly supported continuous function $\dot{v}(x, y)$ on $\R \times (0, \dot{R})$ which is twice continuously differentiable with respect to $x$, such that $\dot{v}$ and $\partial_x \dot{v}$ are weakly differentiable with respect to $y$, and such that $\dot{a}(y) \partial_y \dot{v}(x, y)$ and $\dot{a}(y) \partial_{xy} \dot{v}(x, y)$ are essentially bounded on $\R \times (0, \dot{R})$. These conditions assert that the integrals in~\eqref{eq:harm:3:def} make sense. If $\dot{a}$ and $\dot{b}$ are sufficiently regular (for example, locally bounded on $(0, \dot{R})$), then every smooth, compactly supported function is a suitable test function. We will shortly see that also under our more general assumptions on $\dot{a}$ and $\dot{b}$, the class of suitable test functions for~\eqref{eq:harm:3:def} is sufficiently rich.

Given the parameters $\dot{R}$, $\dot{a}$ and $\dot{b}$ of $\dot{\op}$, we first construct the parameters $R$, $\tilde{a}$ and $\tilde{d}$ of the corresponding equation $\tilde{\op} \tilde{u} = 0$ in the Eckhardt--Kostenko form. This involves an appropriate change of scale. Only later we switch to the standard form $\op u = 0$, with appropriately chosen coefficients $a$ and $b$.

The change of scale is determined by the function
\formula[eq:dotsigma]{
 \dot{\sigma}(y) & = \int_0^y \frac{1}{\dot{a}(r)} \, dr .
}
We define
\formula[eq:var23]{
 R & = \dot{\sigma}(\dot{R}^-) , & \qquad \tilde{a}(\dot{\sigma}(y)) & = (\dot{a}(y))^2 - (\dot{b}(y))^2 , \\
 \tilde{u}(x, \dot{\sigma}(y)) & = \dot{u}(x, y) , & \qquad \tilde{d}(\dot{\sigma}(y)) & = \dot{b}'(y) ,
}
and we suppose that $\dot{u}$ is a harmonic function for $\dot{\op}$. Note that $\tilde{a}(s) \ge 0$ for all $s \in [0, R)$, and that $\tilde{a}(s) ds$ is a locally finite measure on $[0, R)$. As in~\cite{km18}, one shows that the $\leb^2(\R)$-valued function $y \mapsto \tilde{u}(\cdot, y)$ is weakly differentiable on $(0, R)$, and
\formula{
 \partial_y \tilde{u}(\cdot, \dot{\sigma}(y)) & = \dot{a}(y) \partial_y \dot{u}(\cdot, y)
}
for almost all $y \in [0, R)$. Let $v$ be a smooth, compactly supported function on $\R \times (0, R)$, and define $\dot{v}$ by the formula $\dot{v}(x, y) = v(x, \dot{\sigma}(y))$. It is easy to see that $\dot{v}$ is a suitable test function for~\eqref{eq:harm:3:def} (so that, in particular, the class of suitable test functions is rich: it is dense in the space of compactly supported, continuous functions), and, with the above notation, formula~\eqref{eq:harm:3:def} reads
\formula{
 & \int_0^{\dot{R}} \biggl( \int_{-\infty}^\infty \tilde{u}(x, \dot{\sigma}(y)) \partial_{xx} v(x, \dot{\sigma}(y)) dx \biggr) \dot{a}(y) dy \\
 & \qquad + \int_0^{\dot{R}} \biggl( \int_{-\infty}^\infty \tilde{u}(x, \dot{\sigma}(y)) \bigl(\partial_{x y} v(x, \dot{\sigma}(y)) - \dot{b}(y) \partial_{xx} v(x, \dot{\sigma}(y))\bigr) dx \biggr) \frac{\dot{b}(y)}{\dot{a}(y)} \, dy \\
 & \qquad\qquad - \int_0^{\dot{R}} \biggl( \int_{-\infty}^\infty \partial_y \tilde{u}(x, \dot{\sigma}(y)) \bigl(\partial_y v(x, \dot{\sigma}(y)) - \dot{b}(y) \partial_x v(x, \dot{\sigma}(y))\bigr) dx \biggr) \frac{1}{\dot{a}(y)} \, dy = 0 .
}
Substituting $s = \dot{\sigma}(y)$ and noting that $ds = (1 / \dot{a}(y)) dy$, we find that, with $y = \dot{\sigma}^{-1}(s)$,
\formula{
 & \int_0^R \biggl( \int_{-\infty}^\infty \tilde{u}(x, s) \partial_{xx} v(x, s) dx \biggr) (\dot{a}(y))^2 ds \\
 & \qquad + \int_0^R \biggl( \int_{-\infty}^\infty \tilde{u}(x, s) \bigl(\partial_{x y} v(x, s) - \dot{b}(y) \partial_{xx} v(x, s)\bigr) dx \biggr) \dot{b}(y) ds \\
 & \qquad\qquad - \int_0^R \biggl( \int_{-\infty}^\infty \partial_y \tilde{u}(x, s) \bigl(\partial_y v(x, s) - \dot{b}(y) \partial_x v(x, s)\bigr) dx \biggr) ds = 0 .
}
After rearrangement, we eventually obtain
\formula{
 & \int_0^R \biggl( \int_{-\infty}^\infty \tilde{u}(x, s) \partial_{xx} v(x, s) dx \biggr) \tilde{a}(s) ds \\
 & \qquad + \int_0^R \biggl( \int_{-\infty}^\infty \bigl(\partial_y \tilde{u}(x, s) \partial_x v(x, s) + \tilde{u}(x, s) \partial_{x y} v(x, s) \bigr) dx \biggr) \dot{b}(y) ds \\
 & \qquad\qquad - \int_0^R \biggl( \int_{-\infty}^\infty \partial_y \tilde{u}(x, s) \partial_y v(x, s) dx \biggr) ds = 0 ,
}
which is precisely formula~\eqref{eq:harm:2:def} in the definition of a function harmonic with respect to $\tilde{\op}$. Thus, $\tilde{u}$ is a harmonic function for $\tilde{\op}$, in the sense of Definition~\ref{def:harm:2}.

In order to find the corresponding operator $\op$ in standard form, we now use the result obtained earlier in this section. We define the coefficients $a$ and $b$ via the formulae
\formula[eq:var13]{
 a(\dot{\sigma}(y)) & = (\dot{a}(y))^2 , & \qquad b(\dot{\sigma}(y)) & = -\dot{b}(y) ,
}
we let, as usual,
\formula{
 B(\dot{\sigma}(y)) & = \int_0^{\dot{\sigma}(y)} b(s) ds = -\int_0^y \frac{\dot{b}(r)}{\dot{a}(r)} \, dr ,
}
and we define
\formula{
 u(x, \dot{\sigma}(y)) & = \tilde{u}(x - B(\dot{\sigma}(y)), \dot{\sigma}(y)) = \dot{u}(x - B(\dot{\sigma}(y)), y) .
}
Note that the formula $a(dy) = a(y) dy$ defines a locally finite measure on $[0, R)$ with a positive almost everywhere density function $a(y)$. It follows that if $\dot{u}$ is a harmonic function for $\dot{\op}$ in the sense of Definition~\ref{def:harm:3}, then $u$ is a harmonic function for $\op$ in the sense of Definition~\ref{def:harm}.

Suppose now that $\dot{u}$ is a harmonic function for $\dot{\op}$ with boundary values $f(x) = \dot{u}(x, 0)$. According to Definition~\ref{def:dn}, the Dirichlet-to-Neumann operator $\dn$ associated to $\dot{u}$ is given by
\formula[eq:dn:3]{
 \dn f(x) & = \lim_{y \to 0^+} \frac{\dot{u}(x - B(\dot{\sigma}(y)), y) - \dot{u}(x, 0)}{\dot{\sigma}(y)} \\
 & = \lim_{y \to 0^+} \bigl(\dot{a}(y) \partial_y \dot{u}(x, y) + \dot{b}(y) \partial_x \dot{u}(x, y) \bigr) ,
}
with the limits in $\leb^2(\R)$; the second inequality is a consequence of~\eqref{eq:dn:def:alt}. As an immediate corollary of Theorem~\ref{thm:main}, we obtain the following result.

\begin{proposition}
\label{prop:main:3}
\begin{enumerate}[label={\textnormal{(\alph*)}}]
\item Under the assumptions listed in Definition~\ref{def:harm:3}, the Dirichlet-to-Neumann operator $\dn$ associated to the equation $\dot{\op} u = 0$, with $\dot{\op}$ given by~\eqref{eq:op:3}, is an operator of class $\dns$.
\item Every operator $\dn$ of class $\dns$ is the Dirichlet-to-Neumann operator associated to the equation $\dot{\op} \dot{u} = 0$ for at most one triplet of parameters $\dot{R}$, $\dot{a}$ and $\dot{b}$ satisfying the conditions listed in Definition~\ref{def:harm:3}.
\end{enumerate}
\end{proposition}

Note that the counterpart of Theorem~\ref{thm:main}\ref{thm:main:b} is incomplete: not all operators of class $\dns$ can be realised as described above. This is the main reason for us to focus on the equation $\op u = 0$ in standard form studied in Section~\ref{sec:main}. On the other hand, some examples take a particularly simple form when written as in~\eqref{eq:eq:3}, and this form is also more suitable for some constructions; further discussion can be found in Section~\ref{sec:ex}.

It is easy to find an operator $\op$ (or $\tilde{\op}$) which corresponds to a given operator $\dot{\op}$, using~\eqref{eq:dotsigma} and~\eqref{eq:var13} (or~\eqref{eq:dotsigma} and~\eqref{eq:var23}). The converse is slightly more complicated. Let $\op$ be an operator in the standard form~\eqref{eq:op:1}, with coefficients $a$ and $b$, and suppose that $a(dy)$ has a positive almost everywhere density function, denoted by $a(y)$. The coefficients $\dot{a}$ and $\dot{b}$ of the corresponding operator $\dot{\op}$ are clearly given again by~\eqref{eq:var13}, with the function $\dot{\sigma}$ completely determined by~\eqref{eq:dotsigma}. More precisely, formula~\eqref{eq:dotsigma} implies that
\formula{
 \int_0^{\dot{\sigma}(y)} \sqrt{a(s)} ds & = y ,
}
and thus $\dot{\sigma}$ is the inverse function of $y \mapsto \int_0^y \sqrt{a(s)} ds$. It is now easy to check that $\dot{a}$ and $\dot{b}$ satisfy all conditions listed in Definition~\ref{def:harm:3}; we omit the details.

\subsection{Notation and preliminaries}
\label{sec:not}

Throughout the article, all measures are assumed to be locally finite and complex-valued measures. By $\ph(t^+)$ and $\ph(t^-)$ we denote one-sided limits of $\ph$ at $t$. As usual, we denote by $\cont_c^\infty(D)$ the class of smooth, compactly supported functions on $D$, and by $\leb^p(D)$ the class of $p$-integrable Borel functions on $D$, with functions equal almost everywhere identified.

The Fourier transform of a function $f$ is denoted by $\hat{f}$: if $f \in \leb^1(\R)$, then $\hat{f}(\xi) = \int_{-\infty}^\infty e^{-i \xi x} f(x) dx$, and the Fourier transformation $f \mapsto \hat{f}$ is continuously extended to $\leb^2(\R)$. Note that if $g(x) = f(x + a)$ then $\hat{g}(\xi) = e^{i \xi a} \hat{f}(\xi)$, while if $g(x) = f'(x)$, then $\hat{g}(\xi) = i \xi \hat{f}(\xi)$.

If $\ph$ is an absolutely continuous function on an interval, then $\ph$ is differentiable almost everywhere, and the weak (or distributional) derivative of $\ph$ corresponds to a function equal almost everywhere to the point-wise derivative. If $\ph$ is a function of bounded variation, then the distributional derivative of $\ph$ corresponds to a measure. Here we take special care about the endpoints of the domain of $\ph$: if $\ph$ is defined on $[0, R)$ and $\ph(0) \ne \ph(0^+)$, then we understand that $\ph'$ contains an atom at $0$ of mass $\ph(0^+) - \ph(0)$, as if $\ph$ was extended to a constant function $\ph(t) = \ph(0)$ for $t < 0$. In particular, the value of $\ph$ at a single point $0$ does influence the distributional derivative of $\ph$.

If $\ph_1$ and $\ph_2$ are functions of bounded variation with no common discontinuities, then $\ph_1 \ph_2$ is of bounded variation, too, and $(\ph_1 \ph_2)' = \ph_1 \ph_2' + \ph_2 \ph_1'$ (where all derivatives are taken in the sense of distributions, and correspond to appropriate measures).

A locally integrable function $u(x, y)$ is said to be \emph{weakly differentiable} with respect to $x$ if there is a locally integrable function $v(x, y)$ such that
\formula{
 -\int_{-\infty}^\infty \int_{-\infty}^\infty u(x, y) \partial_x w(x, y) dx dy & = -\int_{-\infty}^\infty \int_{-\infty}^\infty v(x, y) w(x, y) dx dy
}
for every smooth, compactly supported (test) function $w$. As remarked above, a function $u$ of one variable is weakly differentiable if and only if it is (locally) absolutely continuous. In higher dimensions, we will use the following characterisation of weak differentiability, known as \emph{absolute continuity on lines} (ACL): $u(x, y)$ is weakly differentiable with respect to $x$ if and only if there is a function $\tilde{u}(x, y)$ which is equal to $u(x, y)$ almost everywhere, which is absolutely continuous with respect to $x$ for every $y$, and such that the point-wise derivative $\partial_x \tilde{u}(x, y)$ (which necessarily exists almost everywhere) is a locally integrable function. In this case $\partial_x \tilde{u}(x, y)$ is the weak derivative of $u(x, y)$.

We use the same notation $\partial_x u$ for both the usual (point-wise) and the weak derivative. Whenever this convention may lead to ambiguities, we will explicitly state which derivative we have in mind.

%
%                            ---------- o ----------
%

\section{Auxiliary ODE}
\label{sec:ode}

As it will become apparent in the next section, Fourier transform reduces our problem to the study of a second-order linear ordinary differential equation
\formula[eq:ode]{
 \ph''(dt) & = \xi^2 a(dt) \ph(t) - 2 i \xi b(t) \ph'(t) dt .
}
Here $\ph$ is a function on $[0, R)$ with $R \in (0, \infty]$, $\xi$ is a `spectral' parameter, and the coefficients $a(dt)$ and $b(t)$ are as in the definition of class $\ops$ (Definition~\ref{def:op}): $a(dt)$ is a non-negative measure on $[0, R)$ (we allow for an atom at $0$), the coefficient $b(t)$ is locally square-integrable on $[0, R)$, and $a(dt) - (b(t))^2 dt$ is assumed to be a non-negative measure on $[0, R)$. For our later needs it is enough to assume that $\xi \in \R$; however, we stress that in the proof of Theorem~\ref{thm:ode} arbitrary complex~$\xi$ need to be considered. We understand~\eqref{eq:ode} in the sense of distributions; more precisely, we assume that $\ph$ is an absolutely continuous function such that the first distributional derivative of $\ph$ corresponds to a left-continuous function (which we denote $\ph'(t)$), the second distributional derivative of $\ph$ is a complex-valued measure (that we denote by $\ph''(dt)$), and we have equality of measures given by~\eqref{eq:ode}.

As already mentioned, for our purposes we only need to study the properties of solutions of~\eqref{eq:ode} when $\xi$ is a real number. It is in fact sufficient to consider $\xi > 0$: if $\ph$ is a solution of~\eqref{eq:ode} for some $\xi > 0$, then $\overline{\ph}$ satisfies~\eqref{eq:ode} with $\xi$ replaced by $-\xi$. Furthermore, for $\xi = 0$, equation~\eqref{eq:ode} requires $\ph$ to be an affine function. For this reason, we restrict our attention to $\xi > 0$ in the following statement. We refer to~\cite{ek16} and to Appendix~\ref{app:proof} for results that cover general complex $\xi$.

The following statement summarizes some of the main results of~\cite{ek16}, which play a crucial role in our development. The function $\sym(\xi) / \xi$ introduced in item~\ref{thm:ode:b} is often called the \emph{principal Weyl--Titchmarsh function} for the equation~\eqref{eq:ode}.

\begin{theorem}
\label{thm:ode}
\begin{enumerate}[label={\textnormal{(\alph*)}}]
\item\label{thm:ode:a} Suppose that the coefficients $a(dt)$ and $b(t)$, defined on $[0, R)$, satisfy the conditions of Definition~\ref{def:op}. For every $\xi > 0$ there is a unique solution $\ph_\xi$ of~\eqref{eq:ode} on $[0, R)$ which satisfies $\ph_\xi(0) = 1$ and such that $\ph_\xi$ is bounded when $R = \infty$ (in this case every other solution diverges to infinity at $\infty$), and $\ph_\xi(R^-) = 0$ if $R < \infty$ (in this case every other solution is bounded away from zero in some left neighbourhood of $R$).
\item\label{thm:ode:b} If $\ph_\xi$ is the solution defined above and $\sym(\xi) = -\ph_\xi'(0)$, then $\sym$ extends to a Rogers function; that is, $\sym$ has a holomorphic extension to the the right complex half-plane, and this extension satisfies $\re(\sym(\xi) / \xi) \ge 0$ whenever $\re \xi > 0$.
\item\label{thm:ode:c} If $\ph_\xi$ is the solution defined above, then $|\ph_\xi|^2$ is positive, non-increasing and convex on $[0, R)$, and $|\ph_\xi'|$ is non-increasing on $[0, R)$; furthermore, if $\tilde{a}(dt) = a(dt) - (b(t))^2 dt$, $B(t) = \int_0^t b(s) ds$ and $\tilde{\ph}_\xi(t) = e^{i \xi B(t)} \ph_\xi(t)$, then for every $t \in [0, R)$ we have
\formula{
 \xi^2 \int_{[t, R)} |\tilde{\ph}_\xi(s)|^2 \tilde{a}(ds) + \int_t^R |\tilde{\ph}_\xi'(s)|^2 ds & \le \min \biggl(\re \sym(\xi), \frac{2}{t}\biggr) .
}
In particular, $|\tilde{\ph}_\xi'|^2$ is integrable on $[0, R)$.
\item\label{thm:ode:d} To every Rogers function $\sym$ there corresponds exactly one pair of coefficients $a(dt)$ and $b(t)$, defined on some $[0, R)$ with $R \in (0, \infty]$.
\end{enumerate}
\end{theorem}

In~\cite{ek16}, Eckhardt and Kostenko study the equation~\eqref{eq:ode} in a different form, for the function $\tilde{\ph}_\xi$ rather than $\ph_\xi$. For this reason, we include below a brief discussion of equivalence of these two forms. The direct part of Theorem~\ref{thm:ode} (that is, items~\ref{thm:ode:a} through~\ref{thm:ode:c}) is proved in Sections~3--5 of~\cite{ek16}. For a less general class of coefficients $a(dt)$ and $b(t)$, this goes back to~\cite{kl79,lw98}. The inverse part of Theorem~\ref{thm:ode} (item~\ref{thm:ode:d}) is the main contribution of~\cite{ek16}; its proof involves deep ideas due to de Branges. For reader's convenience, in Appendix~\ref{app:proof} we include an alternative, less abstract proof of parts~\ref{thm:ode:a} through~\ref{thm:ode:c} of Theorem~\ref{thm:ode}, written in the language of~\eqref{eq:ode} rather than that of~\cite{ek16}.

\begin{proof}[Proof of equivalence of Theorem~\ref{thm:ode} and the results of~\cite{ek16}]
We transform equation~\eqref{eq:ode} in a way that corresponds to shearing in Section~\ref{sec:red}: as usual, we denote $B(t) = \int_0^t b(s) ds$, and whenever $\ph$ is a functon on $[0, R)$, we write
\formula{
 \tilde{\ph}(t) & = e^{i \xi B(t)} \ph(t) .
}
On a formal level, $\ph$ is a solution of~\eqref{eq:ode} if and only if $\tilde{\ph}$ satisfies
\formula{
 \tilde{\ph}'' & = e^{i \xi B} (\ph'' + 2 i \xi b \ph' - \xi^2 b^2 \ph + i \xi b' \ph) \\
 & = e^{i \xi B} (\xi^2 (a - b^2) \ph + i \xi b' \ph) = (\xi^2 \tilde{a} + i \xi \tilde{d}) \tilde{\ph} ,
}
where we have denoted $\tilde{a}(dt) = a(dt) - (b(t))^2 dt$ and $\tilde{d}(t) = b'(t)$. By assumption, $\tilde{a}$ is a non-negative measure on $[0, R)$. However, $b$ is only assumed to be locally square-integrable, and therefore the distributional derivative $b'$ need not correspond to a function or a measure: it is an element of the Sobolev space $H_\loc^{-1}([0, R))$ on $[0, R)$ with negative index $-1$, that is, the dual of the Sobolev space $H_c^1([0, R))$ of compactly supported and weakly differentiable functions $f$ on $[0, R)$ such that $f$ and $f'$ are in $\leb^2([0, R))$.

Under the above assumptions (that is, $\tilde{a}$ a non-negative measure and $\tilde{d}$ an element of the Sobolev space $H_\loc^{-1}([0, R))$), the equation satisfied by $\tilde{\ph}$:
\formula[eq:ek]{
 \tilde{\ph}'' & = (\xi^2 \tilde{a} + i \xi \tilde{d}) \tilde{\ph} ,
}
is precisely the equation studied systematically by Eckhardt and Kostenko in~\cite{ek16}, see equation~(1.2) therein. With the notation used there, $z$, $\nu$ and $\omega$ in~\cite{ek16} correspond to $i \xi$, $\tilde{a}$ and $-\tilde{d}$ used here, respectively.

Equivalence of~\eqref{eq:ode} and~\eqref{eq:ek} can be rigorously proved by writing both equations in an integral form. Indeed, suppose that $\tilde{\ph}$ solves~\eqref{eq:ek}. In~\cite{ek16}, this is understood as
\formula{
 -C g(0) - \int_0^R \tilde{\ph}'(t) g'(t) dt & = \xi^2 \int_{[0, R)} \tilde{\ph}(t) g(t) \tilde{a}(dt) + i \xi \tilde{d}(\tilde{\ph} g)
}
for every test function $g$ in $H_c^1([0, R))$ and some constant $C$; see Definition~3.1 in~\cite{ek16}. Recall that $\tilde{a}(dt) = a(dt) - (b(t))^2 dt$, and $\tilde{d} = b'$, that is, by definition,
\formula{
 \tilde{d}(\tilde{\ph} g) & = -\int_0^R (\tilde{\ph} g)'(t) b(t) dt .
}
It follows that
\formula{
 -C g(0) - \int_0^R \tilde{\ph}'(t) g'(t) dt \\
 & \hspace*{-9em} = \xi^2 \int_{[0, R)} \tilde{\ph}(t) g(t) a(dt) - \xi^2 \int_{[0, R)} \tilde{\ph}(t) g(t) (b(t))^2 dt - i \xi \int_0^R (\tilde{\ph} g)'(t) b(t) dt .
}
We have $\tilde{\ph}(t) = e^{i \xi B(t)} \ph(t)$, and we write $g(t) = e^{-i \xi B(t)} \tilde{g}(t)$. Note that $g \in H_c^1([0, R))$ if and only if $\tilde{g} \in H_c^1([0, R))$. Since $\tilde{\ph}'(t) = e^{i \xi B(t)} (\ph'(t) + i \xi b(t) \ph(t))$ and $g'(t) = e^{-i \xi B(t)} (\tilde{g}'(t) - i \xi b(t) \tilde{g}(t))$, we find that
\formula{
 -C \tilde{g}(0) - \int_0^R (\ph'(t) + i \xi b(t) \ph(t)) (\tilde{g}'(t) - i \xi b(t) \tilde{g}(t)) dt & \\ 
 & \hspace*{-19em} = \xi^2 \int_{[0, R)} \ph(t) \tilde{g}(t) a(dt) - \xi^2 \int_{[0, R)} \ph(t) \tilde{g}(t) (b(t))^2 dt - i \xi \int_0^R (\ph \tilde{g})'(t) b(t) dt .
}
After simplification, we obtain
\formula{
 -C \tilde{g}(0) - \int_0^R \ph'(t) \tilde{g}'(t) dt & = \xi^2 \int_{[0, R)} \ph(t) \tilde{g}(t) a(dt) - 2 i \xi \int_0^R \ph'(t) \tilde{g}(t) b(t) dt
}
for every $\tilde{g} \in H_c^1([0, R))$. By taking $\tilde{g}(s) = -(t - s) \ind_{[0, t)}(s)$, we find that
\formula{
 \int_0^t \ph'(s) ds & = C t + \xi^2 \int_{[0, t)} (t - s) \ph(s) a(ds) - 2 i \xi \int_0^t (t - s) \ph'(s) b(s) ds
}
for every $t \in [0, R)$. Finally, differentiation leads to
\formula{
 \ph'(t) & = C + \xi^2 \int_{[0, t)} \ph(s) a(ds) - 2 i \xi \int_0^t \ph'(s) b(s) ds ,
}
which is clearly equivalent to~\eqref{eq:ode}. By essentially reversing the steps of the above argument, we find that if $\ph$ satisfies~\eqref{eq:ode}, then $\tilde{\ph}$ is a solution of~\eqref{eq:ek} (we omit the details), and it follows that~\eqref{eq:ode} and~\eqref{eq:ek} are indeed equivalent.

Part~\ref{thm:ode:a} of the theorem is now essentially Lemma~4.2 in~\cite{ek16}, part~\ref{thm:ode:b} follows from Lemma~5.1 in~\cite{ek16}, and part~\ref{thm:ode:c} is essentially given in the proof of Lemma~5.1 in~\cite{ek16} (see the last display in p.~954 therein). As mentioned above, alternative proofs are given in Appendix~\ref{app:proof}. Finally, part~\ref{thm:ode:d} is stated as Theorem~6.1 in~\cite{ek16}.
\end{proof}

In the remaining part of the article, we denote by $\ph_\xi$ the solution of~\eqref{eq:ode} described by Theorem~\ref{thm:ode} if $\xi > 0$, a similar solution $\ph_\xi(t) = \overline{\ph_{-\xi}(t)}$ if $\xi < 0$, and the constant solution $\ph_0(t) = 1$ if $\xi = 0$.

%
%                            ---------- o ----------
%

\section{Harmonic extensions}
\label{sec:pde}

In this section we describe the class of functions harmonic with respect to operators $\op$ of class $\ops$ in terms of Fourier transform and solutions $\ph_\xi(t)$ of the ODE~\eqref{eq:ode}, described in Theorem~\ref{thm:ode}.

We assume, as in Definition~\ref{def:op}, that $a(dy)$ is a non-negative measure on $[0, R)$, $b(y)$ is a locally square-integrable real-valued function on $[0, R)$, and $a(dy) - (b(y))^2 dy$ is non-negative. We commonly use the auxiliary measure $\tilde{a}(dy) = a(dy) - (b(y))^2 dy$ and function $B(y) = \int_0^y b(t) dt$.

We study functions $u(x, y)$ on $\R \times [0, R)$ which are harmonic with respect to the elliptic operator $\op$ in the sense of Definition~\ref{def:harm}. We denote by $\hat{u}(\xi, y)$ the Fourier transform of $u(x, y)$ in variable $x$, whenever well-defined. We equally often work with the function $\tilde{u}(x, y) = u(x + B(y), y)$. Observe that $\hat{\tilde{u}}(\xi, y) = e^{i \xi B(y)} \hat{u}(\xi, y)$.

We begin with the proof of Proposition~\ref{prop:harm}, which asserts the existence and uniqueness of harmonic extensions. The argument is divided into two steps, which correspond to uniqueness and existence, respectively.

\begin{lemma}
\label{lem:ext:1}
Suppose that $\op$ is an operator of class $\ops$. For $\xi \in \R$ let $\ph_\xi$ be the solution of~\eqref{eq:ode} discussed in Section~\ref{sec:ode}. If $u$ is harmonic with respect to $\op$, then for all $y \in [0, R)$ we have, for almost all $\xi \in \R$,
\formula{
 \hat{u}(\xi, y) & = \hat{u}(\xi, 0) \ph_\xi(y) .
}
\end{lemma}

\begin{proof}
By Definition~\ref{def:harm} and Plancherel's theorem, $y \mapsto \hat{u}(\cdot, y)$ is again a bounded, continuous mapping from $[0, R)$ to $\leb^2(\R)$, which vanishes at $R^-$ if $R < \infty$. Furthermore, $y \mapsto \hat{\tilde{u}}(\cdot, y)$ is weakly differentiable on $(0, R)$, and $\partial_y \hat{\tilde{u}}(\cdot, y)$ is the Fourier transform of $\partial_y \tilde{u}(\cdot, y)$ for almost all $y \in (0, R)$ (here $\tilde{u}(x, y) = u(x + B(y), y)$). Our goal is to prove that $\hat{u}(\xi, \cdot)$ is a solution to the ODE~\eqref{eq:ode}. The proof is rather straightforward, but it requires some care due to possible irregularities of $\hat{u}$.

Here is the philosophy of the proof: if $u$ is sufficiently regular, then, by~\eqref{eq:harm:def:weak} and Plancherel's theorem, we have
\formula{
 & -\int_{(0, R)} \int_{-\infty}^\infty \xi^2 \hat{u}(\xi, y) \overline{\hat{v}(\xi, y)} d\xi a(dy) \\
 & \qquad + 2 \int_0^R \int_{-\infty}^\infty i \xi \partial_y \hat{u}(\xi, y) \overline{\hat{v}(\xi, y)} b(y) d\xi dy \\
 & \qquad\qquad - \int_0^R \int_{-\infty}^\infty \partial_y \hat{u}(\xi, y) \overline{\partial_y \hat{v}(\xi, y)} d\xi dy = 0
}
for every $v \in \cont_c^\infty(\R \times (0, R))$. By a density argument, this implies that (after a modification on a set of zero Lebesgue measure) for almost all $\xi \in \R$ the function $\hat{u}(\xi, \cdot)$ is a solution of the ODE~\eqref{eq:ode}, and hence $\hat{u}(\xi, y) = \hat{u}(\xi, 0) \ph_\xi(y)$, as desired. Our goal is to make the above idea rigorous in the general case, where only minimal smoothness of $u$ is assumed.

By Definition~\ref{def:harm} (or, more precisely, by~\eqref{eq:harm:def}) and Plancherel's theorem, for every $v \in \cont_c^\infty(\R \times (0, R))$ we have 
\formula[eq:harm:aux1]{
 & -\int_{(0, R)} \int_{-\infty}^\infty \xi^2 \hat{u}(\xi, y) \overline{\hat{v}(\xi, y)} d\xi \tilde{a}(dy) \\
 & \qquad + 2 \int_0^R \int_{-\infty}^\infty i \xi e^{-i \xi B(y)} \partial_y \hat{\tilde{u}}(\xi, y) \overline{\hat{v}(\xi, y)} b(y) d\xi dy \\
 & \qquad\qquad + \int_0^R \int_{-\infty}^\infty \hat{u}(\xi, y) \overline{\partial_{yy} \hat{v}(\xi, y)} d\xi dy = 0 .
}
The the ACL characterisation of weak differentiability implies that, after modifying $\hat{\tilde{u}}(\xi, y)$ and $\partial_y \hat{\tilde{u}}(\xi, y)$ on a set of zero Lebesgue measure, we may assume that for every $\xi \in \R$ the function $\hat{\tilde{u}}(\xi, \cdot)$ is absolutely continuous on $[0, R)$, and the point-wise derivative of this uni-variate function agrees almost everywhere on $[0, R)$ with the weak derivative $\partial_y \hat{\tilde{u}}(\xi, \cdot)$ of the bi-variate function. We temporarily work with this modification, and a similar modification of $\hat{u}(\xi, y) = e^{i \xi B(y)} \hat{\tilde{u}}(\xi, y)$.

For every $\xi \in \R$, the function $\tilde{\ph}(y) = \hat{\tilde{u}}(\xi, y)$ is absolutely continuous on $[0, R)$. It follows that also $\ph(y) = \hat{u}(\xi, y) = e^{-i \xi B(y)} \hat{\tilde{u}}(\xi, y) = e^{-i \xi B(y)} \tilde{\ph}(y)$ is absolutely continuous on $[0, R)$, and $\tilde{\ph}'(y) = e^{i \xi B(y)} (\ph'(y) + i \xi b(y) \ph(y))$ for almost all $y \in [0, R)$ (we stress that $u(x, y)$ need not be weakly differentiable with respect to $y$; nevertheless, it turns out that $\hat{u}(\xi, y)$ is necessarily weakly differentiable with respect to $y$). Applying this identity to~\eqref{eq:harm:aux1}, we find that
\formula{
 & -\int_{(0, R)} \int_{-\infty}^\infty \xi^2 \hat{u}(\xi, y) \overline{\hat{v}(\xi, y)} d\xi a(dy) \\
 & \qquad + 2 \int_0^R \int_{-\infty}^\infty i \xi \partial_y \hat{u}(\xi, y) \overline{\hat{v}(\xi, y)} b(y) d\xi dy \\
 & \qquad\qquad + \int_0^R \int_{-\infty}^\infty \hat{u}(\xi, y) \overline{\partial_{yy} \hat{v}(\xi, y)} d\xi dy = 0 .
}
We choose $v(x, y) = v_1(x) v_2(y)$ with $v_1 \in \cont_c^\infty(\R)$ and $v_2 \in \cont_c^\infty((0, R))$, so that $\hat{v}(\xi, y) = \hat{v}_1(\xi) v_2(y)$. Using Fubini's theorem, we find that
\formula{
 \int_{-\infty}^\infty \overline{\hat{v}_1(\xi)} & \biggl( - \int_{(0, R)} \xi^2 \hat{u}(\xi, y) v_2(y) a(dy) \\
 & \qquad \phantom{\biggl(} + 2 \int_0^R i \xi \partial_y \hat{u}(\xi, y) v_2(y) b(y) dy \\
 & \qquad\qquad \phantom{\biggl(} + \int_0^R \hat{u}(\xi, y) v_2''(y) dy \biggr) d\xi = 0 .
}
The class of Fourier transforms $\hat{v}_1$ of functions $v_1 \in \cont_c^\infty(\R)$ is dense in $\leb^2(\R)$. Therefore, if $v_2 \in \cont_c^\infty((0, R))$, then for almost all $\xi \in \R$ we have
\formula[eq:harm:aux2]{
 & -\int_{(0, R)} \xi^2 \hat{u}(\xi, y) v_2(y) a(dy) \\
 & \qquad + 2 \int_0^R i \xi \partial_y \hat{u}(\xi, y) v_2(y) b(y) dy \\
 & \qquad\qquad + \int_0^R \hat{u}(\xi, y) v_2''(y) dy = 0 .
}
By choosing a countable, dense set of $v_2 \in \cont_c^\infty((0, R))$, we conclude that for almost all $\xi \in \R$, the above equality is satisfied for a dense set of $v_2 \in \cont_c^\infty((0, R))$, and therefore for all $v_2 \in \cont_c^\infty((0, R))$.

For a fixed $\xi \in \R$ with the above property, we let $\ph(y) = \hat{u}(\xi, \cdot)$. Identity~\eqref{eq:harm:aux2} reads
\formula{
 -\int_{(0, R)} \xi^2 \ph(y) v_2(y) a(dy) + 2 \int_0^R i \xi \ph'(y) v_2(y) b(y) dy + \int_0^R \ph(y) v_2''(y) dy = 0
}
for all $v_2 \in \cont_c^\infty((0, R))$, which is the distributional formulation of the ODE
\formula[eq:ode:h]{
 \xi^2 \ph(y) a(dy) - 2 i \xi b(y) \ph'(y) dy - \ph''(dy) & = 0 ,
}
identical to~\eqref{eq:ode}, studied in the previous section.

Suppose that $R = \infty$. By Theorem~\ref{thm:ode}, in this case any solution of~\eqref{eq:ode:h} is either a multiple of $\ph_\xi$ or it diverges to infinity at $R^-$. Since the $\leb^2(\R)$ norm of $\hat{u}(\cdot, y)$ is bounded uniformly with respect to $y$ in $[0, \infty)$ except a set of zero Lebesgue measure (recall that we have modified $\hat{u}$ on a set of zero Lebesgue measure!), $|\hat{u}(\xi, y)|$ cannot diverge to infinity as $y \to \infty$ for all $\xi$ in a set of positive Lebesgue measure (otherwise, by Fatou's lemma, the $\leb^2(\R)$ norm of $\hat{u}(\cdot, y)$ would diverge to infinity as $y \to \infty$). It follows that for almost all $\xi \in \R$ there is a number $c_\xi \in \C$ such that for all $y \in (0, R)$ we have
\formula{
 \hat{u}(\xi, y) & = c_\xi \ph_\xi(y) .
}
The same equality necessarily holds almost everywhere for the original version of $\hat{u}$, before modification on a set of zero Lebesgue measure. Since $y \mapsto \hat{u}(\cdot, y)$ is a continuous map from $[0, R)$ to $\leb^2(\R)$, and $\ph_\xi(0) = 1$ for all $\xi \in \R$, we have $c_\xi = \hat{u}(0, \xi)$ for almost all $\xi \in \R$, and the assertion of the lemma follows.

When $R < \infty$, the proof is very similar. In this case we know that the $\leb^2(\R)$ norm of $\hat{u}(\cdot, y)$ converges to zero as $y \to R^-$ except for a set of $y$ of zero Lebesgue measure, and by Theorem~\ref{thm:ode}, any solution $\ph$ of~\eqref{eq:ode:h} is either a multiple of $\ph_\xi$ or $|\ph|$ has a positive lower limit at $R^-$. Fatou's lemma again implies that for almost every $\xi \in \R$ the function $\hat{u}(\xi, \cdot)$ is a multiple of $\ph_\xi$, and the remaining part of the argument is the same as in the case $R = \infty$.
\end{proof}

\begin{lemma}
\label{lem:ext:2}
Suppose that $\op$ is an operator of class $\ops$. For $\xi \in \R$ let $\ph_\xi$ be the solution of~\eqref{eq:ode} discussed in Section~\ref{sec:ode}. If $f \in \leb^2(\R)$, then the formula
\formula{
 \hat{u}(\xi, y) & = \hat{f}(\xi) \ph_\xi(y)
}
defines a function $u$ on $\R \times [0, R)$ harmonic with respect to $\op$.
\end{lemma}

\begin{proof}
We need to verify the conditions listed in Definition~\ref{def:harm}. By Theorem~\ref{thm:ode}, for every $\xi \in \R$ the function $\ph_\xi$ is continuous and bounded by $1$. In particular, for every $y \in [0, R)$, $\hat{u}(\cdot, y)$ is in $\leb^2(\R)$ with norm bounded by the $\leb^2(\R)$ norm of $\hat{f}$, and so it is the Fourier transform of some function $u(\cdot, y)$ with $\leb^2(\R)$ norm no greater than the $\leb^2(\R)$ norm of $f$. Since $\ph_\xi$ is continuous on $[0, R)$ for every $\xi \in \R$, by the dominated convergence theorem, $y \mapsto \hat{u}(\cdot, y)$ is a continuous map from $[0, R)$ into $\leb^2(\R)$; thus $y \mapsto u(\cdot, y)$ has the same property. A similar argument implies that if $R < \infty$, then $u(\cdot, y)$ converges in $\leb^2(\R)$ to zero as $y \to R^-$. This proves that condition~\ref{def:harm:a} of Definition~\ref{def:harm} is satisfied.

As usual, let $B(y) = \int_0^y b(t) dt$ and $\tilde{u}(x, y) = u(x + B(y), y)$, so that
\formula{
 \hat{\tilde{u}}(\xi, y) & = e^{i \xi B(y)} \hat{u}(\xi, y) = \hat{f}(\xi) \tilde{\ph}_\xi(y) ,
}
where $\tilde{\ph}_\xi(y) = e^{i \xi B(y)} \ph_\xi(y)$. By Theorem~\ref{thm:ode}, for every $\xi \in \R$ and $y \in (0, R)$, the function $\tilde{\ph}_\xi$ is weakly differentiable, and $\tilde{\ph}_\xi'$ is square integrable on $[y, R)$, with $\leb^2((y, R))$ norm bounded by $1 / \sqrt{2 y}$. Therefore, if we define
\formula{
 \partial_y \hat{\tilde{u}}(\xi, y) & = \hat{f}(\xi) \tilde{\ph}_\xi'(y) ,
}
then $\partial_y \hat{\tilde{u}}$ is square integrable on $\R \times [y, R)$ for every $y \in (0, R)$. Fubini's theorem asserts that for every $v \in C_c^\infty(\R \times (0, R))$ we have
\formula{
 \int_0^R \int_{-\infty}^\infty \partial_y \hat{\tilde{u}}(\xi, y) \overline{\hat{v}(\xi, y)} d\xi dy & = -\int_0^R \int_{-\infty}^\infty \hat{\tilde{u}}(\xi, y) \partial_y \overline{\hat{v}(\xi, y)} d\xi dy ,
}
and by Plancherel's theorem we find that formula~\eqref{eq:weak} is satisfied with $\partial_y \tilde{u}(x, y)$ defined as the inverse Fourier transform of $\partial_y \hat{\tilde{u}}(\xi, y)$. We have already observed that $\partial_y \tilde{u}(x, y)$ is square integrable in every strip $\R \times (y_1, y_2)$ with $0 < y_1 < y_2 < R$, and hence condition~\ref{def:harm:b} of Definition~\ref{def:harm} is satisfied.

Finally, condition~\ref{def:harm:c} of Definition~\ref{def:harm} reduces to an application of Plancherel's theorem and Fubini's theorem. Indeed, by Plancherel's theorem, up to a factor $(2 \pi)^{-1}$, the left-hand side of~\eqref{eq:harm:def} is equal to
\formula{
 & -\int_{(0, R)} \biggl(\int_{-\infty}^\infty \xi^2 \hat{u}(\xi, y) \overline{\hat{v}(\xi, y)} d\xi\biggr) \tilde{a}(dy) \\
 & \qquad + 2 \int_0^R \biggl(\int_{-\infty}^\infty i \xi e^{-i \xi B(y)} \partial_y \hat{\tilde{u}}(\xi, y) \overline{\hat{v}(\xi, y)} d\xi\biggr) b(y) dy \\
 & \qquad\qquad + \int_0^R \biggl(\int_{-\infty}^\infty \hat{u}(\xi, y) \partial_{yy} \overline{\hat{v}(\xi, y)} d\xi\biggr) dy ,
}
which, by Fubini's theorem, is equal to
\formula{
 \int_{-\infty}^\infty & \biggl(-\int_{(0, R)} \xi^2 \hat{u}(\xi, y) \overline{\hat{v}(\xi, y)} \tilde{a}(dy) \\
 & \qquad \phantom{\biggl(} + 2 \int_0^R i \xi e^{-i \xi B(y)} \partial_y \hat{\tilde{u}}(\xi, y) \overline{\hat{v}(\xi, y)} b(y) dy \\
 & \qquad\qquad \phantom{\biggl(} + \int_0^R \hat{u}(\xi, y) \partial_{yy} \overline{\hat{v}(\xi, y)} dy \biggr) d\xi .
}
Using the definitions $\hat{u}(\xi, y) = \hat{f}(\xi) \ph_\xi(y)$ and $\hat{\tilde{u}}(\xi, y) = e^{i \xi B(y)} \hat{f}(\xi) \ph_\xi(y)$, together with the fact that $\ph_\xi$ is a solution of~\eqref{eq:ode}, we find that the expression under the outer integral is zero for every $\xi \in \R$, and hence condition~\ref{def:harm:c} of Definition~\ref{def:harm} is satisfied.
\end{proof}

The above two lemmas prove Proposition~\ref{prop:harm}.

\begin{lemma}
\label{lem:dn1}
Suppose that $\op$ is an operator of class $\ops^\star$. For $\xi \in \R$ let $\ph_\xi$ be the solution of~\eqref{eq:ode} discussed in Section~\ref{sec:ode}, and let $\sym(\xi) = -\ph_\xi'(0)$ be the associated Rogers function. If $u$ is a harmonic function for $\op$ (in the sense of Definition~\ref{def:harm}) with boundary values $f \in \leb^2(\R)$, then the $\leb^2(\R)$ limit in the definition of the Dirichlet-to-Neumann operator
\formula[eq:dn:def1]{
 \dn f & = \partial_y u(\cdot, 0) = \lim_{y \to 0^+} \frac{u(\cdot, y) - u(\cdot, 0)}{y}
}
exists if and only if $\sym(\xi) \hat{f}(\xi)$ is square integrable, and in this case
\formula[eq:dn1]{
 \widehat{\dn f}(\xi) & = -\sym(\xi) \hat{f}(\xi) .
}
\end{lemma}

\begin{proof}
By Lemma~\ref{lem:ext:1}, for every $y \ge 0$ and $\xi \in \R$ we have $\hat{u}(\xi, y) = \hat{f}(\xi) \ph_\xi(y)$ (after choosing the right representative of $\hat{u}(\cdot, y)$); and conversely, by Lemma~\ref{lem:ext:2}, for every $f \in \leb^2(\R)$ there is a corresponding function $u$ harmonic with respect to $\op$. By Theorem~\ref{thm:ode}, $|\ph_\xi'|$ is non-increasing, so that $|\ph_\xi'(y)| \le |\ph_\xi'(0)| = |\sym(\xi)|$ for all $y \in [0, R)$. It follows that
\formula[eq:dn:def1:f]{
 \lim_{y \to 0^+} \frac{\hat{u}(\xi, y) - \hat{u}(\xi, 0)}{y} & = \hat{f}(\xi) \ph'_\xi(0) = -\hat{f}(\xi) \sym(\xi)
}
for every $\xi \in \R$, and 
\formula{
 \biggl| \frac{\hat{u}(\xi, y) - u(\xi, 0)}{y} \biggr| & \le |\sym(\xi) \hat{f}(\xi)|
}
for every $\xi \in \R$ and $y \in (0, R)$. If $\sym \hat{f} \in \leb^2(\R)$, then, by dominated convergence, the limit in~\eqref{eq:dn:def1:f} exists in $\leb^2(\R)$. By Plancherel's theorem, the limit in~\eqref{eq:dn:def1} exists in $\leb^2(\R)$, and~\eqref{eq:dn1} holds. Conversely, if the limit in~\eqref{eq:dn:def1} exists in $\leb^2(\R)$, then, again by Plancherel's theorem, the limit in~\eqref{eq:dn:def1:f} exists in $\leb^2(\R)$, and it is necessarily equal to $\sym \hat{f}$. Consequently, $\sym \hat{f} \in \leb^2(\R)$, as desired.
\end{proof}

The above lemma proves the first statement of Theorem~\ref{thm:main} for operators $\op$ of class $\ops^\star$. As explained after the statement of Theorem~\ref{thm:main}, extension to the class $\ops$ is immediate. The other part of Theorem~\ref{thm:main} is a consequence of item~\ref{thm:ode:d} of Theorem~\ref{thm:ode}.

%
%                            ---------- o ----------
%

\section{Examples}
\label{sec:ex}

In this section we discuss a number of non-local operators and corresponding extension problems. More precisely, we prescribe the coefficients $a$ and $b$ of the reduced elliptic equation $\op u = 0$, and evaluate, often omitting the technical details, the corresponding solution $\ph_\xi(t)$ of the ODE~\eqref{eq:ode}. This allows us to identify the corresponding Fourier symbol $-\sym(\xi) = \ph_\xi'(0)$, and eventually leads to the explicit form of the Dirichlet-to-Neumann operator $\dn$. Whenever possible, we discuss all three variants: the standard form $\op$, the Eckhardt--Kostenko form $\tilde{\op}$ and the divergence-like form $\dot{\op}$, discussed in Section~\ref{sec:dis}. For the convenience of the reader, we recall that
\formula{
 \op u & = a(dy) \partial_{xx} u + 2 b(y) \partial_{x y} u + \partial_{yy} u , \\
 \tilde{\op} \tilde{u} & = \tilde{a}(dy) \partial_{xx} \tilde{u} + \partial_{yy} \tilde{u} + \tilde{d}(y) \partial_x \tilde{u} , \\
 \dot{\op} \dot{u} & = \nabla_{x,y} \cdot (\dot{a}(y) \nabla_{x,y} \dot{u}) + 2 \dot{b}(y) \partial_{xy} \dot{u} .
}
We begin with two rather trivial examples, then we discuss three general constructions, and finally we discuss the representation of non-symmetric fractional derivatives.

\subsection{Zero operator}

If $a(dy) = 0 \hspace*{1pt} dy$ and $b(y) = 0$ for all $y \in [0, \infty)$, then the solution of the ODE~\eqref{eq:ode} is given by
\formula{
 \ph_\xi(y) & = 1 ,
}
and consequently
\formula{
 \sym(\xi) & = -\ph_\xi'(0) = 0 & \text{and} && \dn f(x) & = 0 .
}
Therefore, the equation $\op u = 0$ (or $\tilde{\op} \tilde{u} = 0$ with $\tilde{a}(dy) = 0 \hspace*{1pt} dy$ and $b(y) = 0$) in $\R \times [0, \infty)$ corresponds to the Dirichlet-to-Neumann operator $\dn f = 0$.

Note, however, that the same coefficients $a(dy) = 0 \hspace*{1pt} dy$ and $b(y) = 0$ on a finite interval $[0, R)$ lead to a non-zero Dirichlet-to-Neumann operator $\dn$. Indeed, if we set $\gamma = 1 / R$, then we easily find that
\formula{
 \ph_\xi(y) & = 1 - \gamma y , & \qquad \sym(\xi) & = -\ph_\xi'(0) = \gamma , & \qquad \dn f(x) & = -\gamma f(x) .
}
Therefore, the equation $\op u = 0$ (or $\tilde{\op} \tilde{u} = 0$) in $\R \times [0, 1 / \gamma)$ corresponds to the Dirichlet-to-Neumann operator $\dn f = -\gamma f$.

\subsection{Constant coefficients}
\label{sec:const}

Let $p \ge 0$, $q \in \R$, and consider $a(dy) = (p^2 + q^2) dy$ and $b(y) = -q$ for $y \in [0, \infty)$. Then
\formula{
 \ph_\xi(y) & = e^{(-p |\xi| + i q \xi) y} , & \qquad \sym(\xi) & = -\ph_\xi'(0) = p |\xi| - i q \xi .
}
Thus, unsurprisingly, $\op u = (p^2 + q^2) \partial_{xx} u - 2 q \partial_{xy} u + \partial_{yy} u$ corresponds to the Dirichlet-to-Neumann operator
\formula{
 \dn f(x) & = -p (-\partial_{xx})^{1/2} f(x) + q f'(x) .
}
Here $\partial_{xx}$ is the second derivative operator (the one-dimensional Laplace operator), and $(-\partial_{xx})^{1/2}$ is the usual Dirichlet-to-Neumann operator for the Laplace equation in the half-plane. In other words, $\dn$ corresponds to $\alpha = 0$, $\beta = q$, $\gamma = 0$ and $\nu(z) = p \pi^{-1} |z|^{-2}$ in Definition~\ref{def:cm}.

The corresponding operator $\tilde{\op}$ in Eckhardt--Kostenko form is simply $\tilde{\op} \tilde{u} = p^2 \partial_{xx} \tilde{u} + \partial_{yy} \tilde{u}$, with coefficients $\tilde{a}(dy) = p^2 dy$ and $\tilde{d}(y) = 0$ that do not depend on $q$. The first-order term $q f'$ in the expression for $\dn f(x)$ comes from the somewhat artificial definition~\eqref{eq:dn:2} of the Dirichlet-to-Neumann operator: the function $B$ is defined by $\tilde{d}$ up to a linear term only, and we choose $B(y) = -q y$ in order that $B'(y) = -q = b(y)$.

%\subsection{Symmetric operators with additional first-order term}

%As explained in the introduction, if we set $b(y) = 0$, then there is a one-to-one correspondence between coefficients $a(dy)$ and \emph{complete Bernstein functions} $\psi$; namely, the Dirichlet-to-Neumann operator associated to the equation $\op u = 0$ (with coefficients $a(dy)$ and $b(y) = 0$) is $\dn = -\psi(-\Delta)$, with corresponding symbol $-\sym(\xi) = -\psi(\xi^2)$ (see~\cite{km18} for a detailed discussion).

%Suppose that the coefficients $a(dy)$ and $b(y) = 0$ correspond to a complete Bernstein function $\psi$ in the above sense, and fix $q \in \R$. If we set $a_q(dy) = a(dy) + q^2 dy$ and $b_q(y) = q$, then the solution of the ODE~\eqref{eq:ode} with coefficients $a_q$ and $b_q$ is given by $\ph_{q,\xi}(y) = e^{-i q \xi y} \ph_\xi(y)$, where $\ph_\xi$ is the solution of~\eqref{eq:ode} with the above $a(dy)$ and $b(y) = 0$. Consequently,
%\formula{
% \sym_q(\xi) & = \sym(\xi) + i q \xi = \psi(\xi^2) + i q \xi .
%}
%It follows that the Dirichlet-to-Neumann operator corresponding to $\op_q u = (a(dy) + q^2) \partial_{xx} u - 2 q \partial_{xy} u + \partial_{yy} u$ is given by
%\formula{
% \dn_q f(x) & = -\psi(-\Delta) f(x) - q f'(x) .
%}
%The corresponding operator $\tilde{\op}_q$ is simply $\tilde{\op}_q u = a(dy) \partial_{xx} \tilde{u} + \partial_{yy} u$, and it does not depend on $q$. Again, the first-order term $q f'$ is present due to the non-standard definition~\eqref{eq:dn:2} of the Dirichlet-to-Neumann operator.

\subsection{Degenerate equations corresponding to one-sided operators without first-order term}
\label{sec:subord}

As explained in the introduction, there is a one-to-one correspondence between measures $a_0(dy)$ on $[0, R)$ and \emph{complete Bernstein functions} $\psi$. Namely, the Dirichlet-to-Neumann operator associated to the equation $\op u = 0$ with coefficients $a(dy) = a_0(dy)$ and $b(y) = 0$ is $\dn = -\psi(-\partial_{xx})$. By this we mean that the corresponding symbol is equal to $-\sym(\xi) = -\psi(\xi^2)$. We refer to~\cite{km18} for a detailed discussion.

It is known that $\partial_{xx}$ can be replaced by a more general non-positive definite operator $D_x$ acting in variable $x$: every operator of the form $-\psi(-D_x)$ arises as the Dirichlet-to-Neumann map for the equation $a_0(dy) D_x u + \partial_{yy} u = 0$. In particular, we can set $D = -\partial_x$. We refer to~\cite{gms13} for a related discussion.

The above observation indicates that the operator $\dn f = -\psi(\partial_x) f$, corresponding to the symbol $-\sym(\xi) = -\psi(i \xi)$, is the Dirichlet-to-Neumann operator associated with the equation $\tilde{\op} \tilde{u} = 0$, where $\tilde{\op} \tilde{u} = \partial_{yy} \tilde{u} - a_0(dy) \partial_x \tilde{u}$. Note that here it is more convenient to work with the operator $\tilde{\op}$ in the Eckhardt--Kostenko form~\eqref{eq:op:2}, with coefficients
\formula{
 \tilde{a}(dy) & = 0 \hspace*{1pt} dy , & \qquad \tilde{d}(y) = -a_0(dy) .
}
If we denote $A_0(y) = a_0([0, y))$, then the corresponding operator $\op$ in the standard form~\eqref{eq:op} is easily found to have coefficients
\formula{
 a(dy) & = (A_0(y))^2 dy , & \qquad b(y) & = A_0(y) .
}
In a similar way, we can find the corresponding operator $\dot{\op}$ in the divergence-like form~\eqref{eq:op:3}, as long as $A_0(y) = a_0([0, y))$ is strictly positive for $y > 0$. Let
\formula{
 B(y) & = \int_0^y b(s) ds = \int_0^y A_0(s) ds = \int_{[0, y)} (y - s) a_0(ds) .
}
Then $B(y) = \int_0^y \sqrt{a(s)} ds$, so that $\dot{\sigma}(y) = B^{-1}(y)$ (see Section~\ref{sec:dis} for the notation), and consequently
\formula{
 \dot{R} & = B(R^-) , & \qquad \dot{a}(y) & = -\dot{b}(y) = b(B^{-1}(y))
}
for $y \in [0, R)$. In other words,
\formula{
 \dot{R} & = \int_{[0, R)} (R - s) a_0(ds) , & \qquad \dot{a}\biggl(\int_{[0, y)} (y - s) a_0(ds)\biggr) & = a_0([0, y))
}
for $y \in [0, R)$, and $\dot{b}(y) = -\dot{a}(y)$ for $y \in [0, \dot{R})$.

It is not difficult to verify that the Dirichlet-to-Neumann operator $\dn$ associated to the equation $\op u = 0$ (with $\op$ as above) is indeed the operator $-\psi(\partial_x)$. As usual, let $B(y) = \int_0^y b(s) ds$, and let $\ph_{0,\xi}$ be the solution of the ODE~\eqref{eq:ode} with coefficients $a_0(dy)$ and $b_0(y) = 0$, for an arbitrary complex parameter $\xi$. Then, for $\xi > 0$, the formula
\formula{
 \tilde{\ph}_\xi(y) & = \ph_{0,\sqrt{i \xi}}(y)
}
defines a solution of the ODE $\tilde{\ph}'' = (i \xi) a_0(dy) \tilde{\ph}$, and thus (by equivalence of~\eqref{eq:ode} and~\eqref{eq:ek})
\formula{
 \ph_\xi(y) & = e^{-i \xi B(y)} \tilde{\ph}_\xi(y) = e^{-i \xi B(y)} \ph_{0,\sqrt{i \xi}}(y)
}
is a solution of~\eqref{eq:ode} with coefficients $a(dy)$ and $b(y)$ defined above. It is more complicated to show that this is the solution discussed in Section~\ref{sec:ode}, that is, that $\ph_\xi$ is bounded if $R = \infty$, and $\ph_\xi$ has a zero left limit at $R$ when $R < \infty$; we omit the details. Since $B'(0) = b(0) = 0$, we find that the symbol $-\sym(\xi)$ of the corresponding Dirichlet-to-Neuman operator is given by
\formula{
 \sym(\xi) & = -\ph_\xi'(0) = -\ph_{0,\sqrt{i \xi}}'(0) = \psi(i \xi) ,
}
and consequently
\formula{
 \dn f(x) & = -\psi(\partial_x) f(x) .
}
With the notation of Definition~\ref{def:cm}, this operator corresponds to $\alpha = 0$, $\gamma \ge 0$, $\nu$ such that $\nu(z) = 0$ for all $z < 0$, and $\beta \ge \int_0^1 z \nu(z) dz$. In other words,
\formula{
 \dn f(x) & = \check{\beta} f'(x) - \gamma f(x) + \int_0^\infty (f(x + z) - f(x)) \nu(z) dz ,
}
where $\check{\beta}, \gamma \ge 0$, $\nu$ is a completely monotone function on $(0, \infty)$, and $\min\{1, z\} \nu(z)$ is integrable over $(0, \infty)$. The operator $-\dn = \psi(\partial_x)$ can be though of as a (right) generalised fractional derivative of order between $0$ and $1$.

We remark that the condition $a_0([0, y)) > 0$ for every $y > 0$ (required in order to properly define the operator $\dot{\op}$ in divergence-like form) is equivalent to $\psi$ being unbounded on $(0, \infty)$. This follows, for example, from formula~(2.14) in~\cite{kw82}; we omit the details.

\subsection{Complementary equations and operators}
\label{sec:comp}

Following Section~5.7 in~\cite{km18}, where symmetric operators are studied, we say that the operators $\dn$ and $\dn^\sharp$ of class $\dns$ are \emph{complementary}, if their composition $\dn \dn^\sharp$ is equal to $-\partial_{xx}$, the one-dimensional Laplace operator. In terms of the corresponding symbols $-\sym$ and $-\sym^\sharp$, we require that $\sym(\xi) \sym^\sharp(\xi) = \xi^2$ for all $\xi \in \R$. We note that if $\sym$ is a Rogers function, then the formula $\sym^\sharp(\xi) = \xi^2 / \sym(\xi)$ also defines a Rogers function (see Proposition~\ref{prop:rogers}); therefore, every operator $\dn$ of class $\dns$ has a unique complementary operator $\dn^\sharp$ of class $\dns$.

In this part it is convenient to work with the equation in a divergence-like form $\dot{\op} \dot{u} = 0$, where $\dot{\op}$ is given by~\eqref{eq:op:3}. Below we argue that if $\dn$ is the corresponding Dirichlet-to-Neumann operator and $\dn^\sharp$ is an operator complementary to $\dn$, then $\dn^\sharp$ is the Dirichlet-to-Neumann operator associated to the \emph{complementary equation} $\dot{\op}^\sharp \dot{u}^\sharp = 0$, with coefficients
\formula{
 \dot{a}^\sharp(y) & = 1 / \dot{a}(y) , & \dot{b}^\sharp(y) & = -\frac{\dot{b}(y)}{(\dot{a}(y))^2} \, ,
}
The proof of this claim consists of two steps.

First, we observe that if $\dot{u}$ is a harmonic function for $\dot{\op}$, then
\formula{
 \dot{u}^\sharp(x, y) & = \dot{a}(y) \partial_y \dot{u}(x, y) + \dot{b}(y) \partial_x \dot{u}(x, y)
}
is a harmonic function for $\dot{\op}^\sharp$. If the coefficients are smooth, this is almost immediately verified using the expression~\eqref{eq:op:3} for $\dot{\op}$, because the operator $\dot{a}(y) \partial_y + \dot{b}(y) \partial_x$ commutes with $\partial_{xx}$. A rigorous proof in the general case is more involved, and we omit the details.

In the second step, we evaluate the Dirichlet-to-Neumann operator $\dn^\sharp$ associated to the equation $\dot{\op}^\sharp \dot{u} = 0$. We already know that $\dn^\sharp$ is an operator of class $\dns$. Let $f$ be a smooth, compactly supported function, let $\dot{u}$ be the harmonic function for $\dot{\op}$ with boundary values $f$, and let $\dot{u}^\sharp$ be defined as above. Then, by~\eqref{eq:dn:3},
\formula{
 \dot{u}^\sharp(x, 0) & = \lim_{y \to 0^+} \dot{u}^\sharp(x, y) = \lim_{y \to 0^+} (\dot{a}(y) \partial_y + \dot{b}(y) \partial_x) \dot{u}(x, y) = \dn f(x)
}
(with all limits understood in the sense of $\leb^2(\R)$). Therefore, $\dot{u}^\sharp$ is a harmonic extension of $\dn f$ for $\dot{\op}^\sharp$. Again using~\eqref{eq:dn:3}, we find that
\formula{
 \dn^\sharp \dn f(x) & = \lim_{y \to 0^+} \biggl(\frac{1}{\dot{a}(y)} \partial_y - \frac{\dot{b}(y)}{(\dot{a}(y))^2} \partial_x\biggr) \dot{u}^\sharp(x, y)
}
(with the limit again understood in the sense of $\leb^2(\R)$). Using the definitions of $\dot{u}^\sharp$ and $\dot{\op}$, we conclude that
\formula{
 \dn^\sharp \dn f(x) & = \lim_{y \to 0^+} \biggl(\frac{1}{\dot{a}(y)} \partial_y - \frac{\dot{b}(y)}{(\dot{a}(y))^2} \partial_x\biggr) (\dot{a}(y) \partial_y + \dot{b}(y) \partial_x) \dot{u}(x, y) \\
 & = \lim_{y \to 0^+} (\dot{\op}^\sharp - \partial_{xx}) \dot{u}(x, y) = -\partial_{xx} f(x) ,
}
as desired (once again with all limits in $\leb^2(\R)$). As in the first step, we omit the technical details related to regularity of $u$ and $u^\sharp$.

It is instructive to evaluate the corresponding coefficients $a$, $b$, $a^\sharp$ and $b^\sharp$ of the complementary equations $\op u = 0$ and $\op^\sharp u = 0$ in standard form. If
\formula{
 \dot{\sigma}(y) & = \int_0^y \frac{1}{\dot{a}(s)} \, ds , & \qquad \dot{\sigma}^\sharp(y) & = \int_0^y \frac{1}{\dot{a}^\sharp(s)} \, ds = \int_0^y \dot{a}(s) ds ,
}
then
\formula{
 R & = \dot{\sigma}(\dot{R}^-) = \int_0^{\dot{R}} \frac{1}{\dot{a}(y)} \, dy , & \qquad R^\sharp & = \dot{\sigma}^\sharp(\dot{R}^-) = \int_0^{\dot{R}} \dot{a}(y) dy ,
}
and
\formula{
 a(\dot{\sigma}(y)) & = (\dot{a}(y))^2 , & \qquad a^\sharp(\dot{\sigma}^\sharp(y)) & = (\dot{a}^\sharp(y))^2 = \frac{1}{(\dot{a}(y))^2} \, , \\
 b(\dot{\sigma}(y)) & = -\dot{b}(y), & \qquad b^\sharp(\dot{\sigma}^\sharp(y)) & = -\dot{b}^\sharp(y) = \frac{\dot{b}(y)}{(\dot{a}(y))^2} .
}
Here we understand that $a(dy) = a(y) dy$ and $a^\sharp(dy) = a^\sharp(y) dy$.

Note that the functions $\dot{\sigma}$ and $\dot{\sigma}^\sharp$ (describing appropriate change of variables) and coefficients $a$ and $a^\sharp$ only depend on the `symmetric' coefficient $\dot{a}$, and not on the `non-symmetric' coefficient $\dot{b}$. Therefore, just as it was the case for symmetric operators (see Section~5.7 in~\cite{km18}), we have
\formula{
 a([0, \dot{\sigma}(y))) & = \int_0^{\dot{\sigma}(y)} a(s) ds = \int_0^y a(\dot{\sigma}(r)) (\dot{\sigma})'(r) dr = \int_0^y \dot{a}(r) dr = \dot{\sigma}^\sharp(y) , \\
 a^\sharp([0, \dot{\sigma}^\sharp(y))) & = \int_0^{\dot{\sigma}^\sharp(y)} a^\sharp(s) ds = \int_0^y a^\sharp(\dot{\sigma}^\sharp(r)) (\dot{\sigma}^\sharp)'(r) dr = \int_0^y \frac{1}{\dot{a}(r)} \, dr = \dot{\sigma}(y) ,
}
so that $y \mapsto a([0, y))$ and $y \mapsto a^\sharp([0, y))$ is a pair of inverse functions. With the terminology of Krein's spectral theory of strings, this means that $a(dy)$ and $a^\sharp(dy)$ are a pair of \emph{dual strings}.

The above argument only covers a limited class of operators $\dn$, namely those operators which are Dirichlet-to-Neumann maps for equations in the divergence-like form $\dot{\op} \dot{u} = 0$. However, the corresponding result in the standard form~\eqref{eq:op} (involving dual Krein's strings) is fully general. A detailed proof is based on the theory of dual Krein's strings and it falls beyond the scope of the present article.

\subsection{Degenerate equations corresponding to one-sided operators with first-order term}

By combining the results of the previous two subsections, we obtain a representation of generalised (left) fractional derivatives of orders between $1$ and $2$. These operators correspond to symbols $-\sym(\xi) = -\xi^2 / \psi(i \xi) = i \xi \psi^\sharp(i \xi)$, where $\psi$ and $\psi^\sharp$ are complete Bernstein functions satisfying $\psi(\xi) \psi^\sharp(\xi) = \xi$. In other words, we formally have $\op = \partial_x \psi^\sharp(\partial_x)$.

Let $a_0(dy)$ be the coefficient associated to $\psi$, and let $\dot{a}$ and $\dot{b}$ be the coefficients of the equation $\dot{\op} \dot{u} = 0$ associated to $\psi(\partial_x)$, as in Section~\ref{sec:subord}. According to Section~\ref{sec:comp}, the complementary equation $\dot{\op}^\sharp \dot{u} = 0$ has coefficients
\formula{
 \dot{a}^\sharp(y) & = \dot{b}^\sharp(y) = \frac{1}{\dot{a}(y)} = \frac{1}{b(B^{-1}(y))} \, ,
}
where $b(y) = a_0([0, y))$ and $B(y) = \int_0^y b(s) ds = \int_{[0, y)} (y - s) a_0(ds)$. In the previous section we have seen that the associated Dirichlet-to-Neumann operator $\dn^\sharp$ has symbol $-\sym^\sharp(\xi) = -\xi^2 / \psi(i \xi) = i \xi \psi^\sharp(i \xi)$, as desired.

We remark that with the notation of Definition~\ref{def:cm}, the operator $\dn^\sharp$ corresponds to $\alpha^\sharp \ge 0$, $\gamma^\sharp = 0$, $\nu^\sharp$ such that $\nu^\sharp(z) = 0$ for all $z < 0$, and $\beta^\sharp \le \int_1^\infty z \nu(z) dz$; that is,
\formula{
 \dn f(x) & = \alpha^\sharp f''(x) - \check{\beta}^\sharp f'(x) + \int_0^\infty (f(x + z) - f(x) - z f'(x)) \nu^\sharp(z) dz ,
}
where $\alpha^\sharp, \check{\beta}^\sharp \ge 0$, $\nu^\sharp$ is a completely monotone function on $(0, \infty)$, and $\min\{z, z^2\} \nu^\sharp(z)$ is integrable over $(0, \infty)$. A detailed discussion of this construction would take us too far from the main scope of this article, and thus we omit the details.

\subsection{Fractional Laplace operator and non-symmetric fractional derivatives}

As discussed in the introduction, if $\mu \in (0, 2)$, $R = \infty$, $a(y) = C_\mu y^{2/\mu - 2}$ for an appropriate $C_\mu$ and $b(y) = 0$, then the corresponding Dirichlet-to-Neumann operator $\dn$ is the fractional Laplace operator $\dn = -(-\partial_xx)^{\mu/2}$; this is the Caffarelli--Silvestre extension technique; see~\cite{cs07}. A similar representation for one-sided fractional derivatives of order $\mu \in (0, 1)$ was studied in detail in~\cite{bmst16}. Here we extend these results to arbitrary (two-sided, non-symmetric) fractional derivatives of order $\mu$.

We first discuss the standard form~\eqref{eq:op}. Let $\mu \in (0, 2)$ be fixed, and suppose that $R = \infty$,
\formula{
 a(y) & = (p^2 + q^2) y^{2/\mu - 2} , & \qquad b(y) & = -q y^{1/\mu - 1}
}
for some $p \ge 0$ and $q \in \R$. The ODE~\eqref{eq:ode} takes form
\formula[eq:ode:st]{
 \ph_\xi''(y) & = (p^2 + q^2) \xi^2 y^{2/\mu - 2} \ph_\xi(y) + 2 i q \xi y^{1/\mu - 1} \ph_\xi'(y) .
}
Our goal is to show that the corresponding symbol is $\sym(\xi) = (A + i B \sign \xi) |\xi|^\mu$ for $\xi \in \R$, where $A \ge 0$ and $B \in \R$ are constants to be determined. Recall that the symbol $\sym(\xi) = -\ph_\xi'(0)$ satisfies $\sym(-\xi) = \overline{\sym(\xi)}$. Thus, with no loss of generality, we assume that $\xi > 0$, and we will show that $\sym(\xi) = (A + i B) \xi^\mu$.

We first consider $\mu \ne 1$ and $p > 0$, and we write
\formula{
 w & = \frac{(1 - \mu) (p - i q)}{2 p} \, , & z(y) & = 2 \mu p \xi y^{1 / \mu} .
}
In this case, with some effort, one verifies that the solution of~\eqref{eq:ode:st} is given by
\formula{
 \ph_\xi(y) & = \frac{\Gamma(\mu + w)}{\Gamma(\mu)} \, \exp(-(1 - \mu)^{-1} w z(y)) \, U(w , 1 - \mu , z(y)) ,
}
where $U$ denotes the confluent hypergeometric function of the second kind (often denoted by $\Psi$; see Section~9.21 in~\cite{gr07} and Section~6.5 in~\cite{emot53}). Using the asymptotic expansion
\formula{
 U(w, 1 - \mu, z) & = \frac{\Gamma(\mu)}{\Gamma(\mu + w)} + \frac{w \, \Gamma(\mu)}{(1 - \mu) \Gamma(\mu + w)} \, z + \frac{\Gamma(-\mu)}{\Gamma(w)} \, z^\mu + O(z^{\min\{\mu + 1, 2\}})
}
as $z \to 0^+$ (see formulae~9.210.1--2 in~\cite{gr07}), we find that
\formula{
 \ph_\xi(y) & = \biggl(1 - \frac{w z(y)}{1 - \mu}\biggr) \biggl(1 + \frac{w z(y)}{1 - \mu} + \frac{\Gamma(-\mu) \Gamma(\mu + w) (z(y))^\mu}{\Gamma(\mu) \Gamma(w)}\biggr) + O((z(y))^{\min\{\mu + 1, 2\}}) \\
 & = 1 + \frac{\Gamma(-\mu) \Gamma(\mu + w)}{\Gamma(\mu) \Gamma(w)} \, (2 \mu p \xi)^\mu y + O(y^{\min\{1 + 1/\mu, 2 / \mu\}}) .
}
Thus, indeed $\ph_\xi(0) = 1$, and
\formula{
 \sym(\xi) & = -\ph_\xi'(0) = \frac{-\Gamma(-\mu) \Gamma(\mu + w)}{\Gamma(\mu) \Gamma(w)} \, (2 \mu p \xi)^\mu .
}
Using the definition of $w$, we eventually find that
\formula{
 \sym(\xi) & = -\frac{\Gamma(-\mu) \Gamma(\mu + \tfrac{(1 - \mu) (p - i q)}{2 p})}{\Gamma(\mu) \Gamma(\tfrac{(1 - \mu) (p - i q)}{2 p})} \, (2 \mu p)^\mu \xi^\mu
}
for $\xi > 0$.

We now move to the case $\mu \ne 1$, $p = 0$ and $q \ne 0$. Let us denote
\formula{
 \thet & = -\tfrac{\pi}{2} \sign (q (1 - \mu)) , & z(y) & = |q \mu (1 - \mu)| \xi y^{1 / \mu} .
}
By a direct calculation, it can be checked that the solution of~\eqref{eq:ode:st} is equal to
\formula{
 \ph_\xi(y) & = \frac{2}{\Gamma(\mu)} \, e^{i \mu \thet/2} (z(y))^{\mu/2} \exp(-e^{i \thet} (1 - \mu)^{-1} z(y)) K_\mu(2 e^{i \thet/2} (z(y))^{1/2}) ,
}
where $K_\mu$ is the modified Bessel function of the second kind (see Section~8.43 in~\cite{gr07} and Section~6.9.1 in~\cite{emot53}). The asymptotic expansion
\formula{
 2 z^{\mu/2} K_\mu(2 z^{1/2}) & = \Gamma(\mu) + \frac{\Gamma(\mu)}{1 - \mu} \, z + \Gamma(-\mu) z^\mu + O(|z|^{\min\{\mu + 1, 2\}})
}
as $|z| \to 0^+$, $|\arg z| \le \tfrac{\pi}{2}$ (see formulae~8.445 and~8.485 in~\cite{gr07}), leads to
\formula{
 \ph_\xi(y) & = \biggl( 1 - \frac{e^{i \thet} z(y)}{1 - \mu} \biggr) \biggl(1 + \frac{e^{i \thet} z(y)}{1 - \mu} + \frac{\Gamma(-\mu) e^{i \mu t} (z(y))^\mu}{\Gamma(\mu)} \biggr) + O((z(y))^{\min\{\mu + 1, 2\}}) \\
 & = 1 + \frac{\Gamma(-\mu)}{\Gamma(\mu)} \, e^{i \mu \thet} (|q \mu (1 - \mu)| \xi)^\mu y + O(y^{\min\{1 + 1/\mu, 2 / \mu\}}) .
}
Again, we find that $\ph_\xi(0) = 1$, and
\formula{
 \sym(\xi) & = -\ph_\xi'(0) = -\frac{\Gamma(-\mu)}{\Gamma(\mu)} \, e^{i \mu \thet} (|q \mu (1 - \mu)| \xi)^\mu .
}
We conclude that, for arbitrary $q \ne 0$,
\formula{
 \sym(\xi) & = -\frac{\Gamma(-\mu)}{\Gamma(\mu)} \, e^{-(i \pi \mu/2) \sign (q (1 - \mu))} |q \mu (1 - \mu)|^\mu \xi^\mu
}
for $\xi > 0$.

Finally, the case $\mu = 1$ was already dealt with in Section~\ref{sec:const}. In this case we simply have
\formula{
 \ph_\xi(y) & = e^{-(p - i q) \xi y} ,
}
so that
\formula{
 \sym(\xi) & = -\ph_\xi'(0) = (p - i q) \xi
}
for $\xi > 0$.

In each case we have $\sym(\xi) = (A + i B) \xi^\mu$ for $\xi > 0$, for some constants $A \ge 0$ and $B \in \R$, and consequently
\formula{
 \sym(\xi) & = (A + i B \sign \xi) |\xi|^\mu
}
for $\xi \in \R$. When $\mu = 1$, we obtain $\sym(\xi) = A |\xi| + i B \xi$, which easily leads to
\formula{
 \dn f(x) & = \frac{A}{\pi} \int_{-\infty}^\infty \frac{f(x + z) - f(x) - z f'(x) \ind_{(-1, 1)}(z)}{|z|^2} \, dy - B f'(x) .
}
If $\mu \ne 1$, we have
\formula{
 \sym(\xi) & = C_+ (-i \xi)^\mu + C_- (i \xi)^\mu ,
}
where both powers are understood as principal branches, and $C_+ e^{-i \mu \pi/2} + C_- e^{i \mu \pi/2} = A + B i$, that is,
\formula{
 C_+ & = \frac{A}{2 \cos \tfrac{\mu \pi}{2}} - \frac{B}{2 \sin \tfrac{\mu \pi}{2}} \, , & \qquad C_- & = \frac{A}{2 \cos \tfrac{\mu \pi}{2}} + \frac{B}{2 \sin \tfrac{\mu \pi}{2}} \, .
}
If $\mu < 1$, it follows that
\formula{
 \dn f(x) & = -\frac{1}{\Gamma(-\mu)} \int_{-\infty}^\infty \frac{f(x + z) - f(x)}{|z|^{1 + \mu}} \, (C_+ \ind_{(0, \infty)}(z) + C_- \ind_{(-\infty, 0)}(z)) dz ,
}
while for $\mu > 1$ we find that
\formula{
 \dn f(x) & = \frac{1}{\Gamma(-\mu)} \int_{-\infty}^\infty \frac{f(x + z) - f(x) - z f'(x)}{|z|^{1 + \mu}} \, (-C_+ \ind_{(0, \infty)}(z) - C_- \ind_{(-\infty, 0)}(z)) dz ;
}
see, for example, Section~7.1 in~\cite{skm93} or Section~31 in~\cite{s99}.

It is immediate to see that the coefficients $a$ and $b$ of the equation $\op u = 0$ in standard form~\eqref{eq:op} correspond to the coefficients
\formula{
 \tilde{a}(y) & = p^2 y^{2/\mu - 2} , & \qquad \tilde{d}(y) & = \frac{q (1 - \mu)}{\mu} \, y^{1/\mu - 2}
}
of the equation $\tilde{\op} \tilde{u} = 0$ in Eckhardt--Kostenko form~\eqref{eq:op:2}. This leads to a certain simplification of the above expressions for $A$ and $B$ (see below). Similarly, one easily finds that the coefficients of the equation $\dot{\op} \dot{u} = 0$ in the divergence-like form are given by
\formula{
 \dot{a}(y) & = \mu^{\mu - 1} (p^2 + q^2)^{\mu/2} y^{1 - \mu} , & \qquad \dot{b}(y) & = \mu^{\mu - 1} q (p^2 + q^2)^{-1/2 + \mu/2} y^{1 - \mu} ;
}
indeed, $\dot{\sigma}$ is the inverse function of $y \mapsto \int_0^y \sqrt{a}(s) ds = \mu \sqrt{p^2 + q^2} y^{1 / \mu}$, that is, $\dot{\sigma}(y) = \mu^{-\mu} (p^2 + q^2)^{-\mu/2} y^\mu$.

\smallskip

The results of this section can be summarised as follows, with a slightly changed notation: we replace $C_\pm$ by $|C_\pm|$. Dirichlet-to-Neumann operators related to the following elliptic equations:
\formula{
 \begin{gathered}
  (p^2 + q^2) y^{2/\mu - 2} \partial_{xx} u - 2 q y^{1/\mu - 1} \partial_{xy} u + \partial_{yy} u = 0 , \\
  \tilde{p}^2 y^{2/\mu - 2} \partial_{xx} \tilde{u} + \partial_{yy} \tilde{u} + \tilde{q} y^{1/\mu - 2} \partial_x \tilde{u} = 0 , \\
  \partial_{xx} \dot{u} + y^{\mu - 1} (\dot{p} \partial_y - \dot{q} \partial_y) \bigl(y^{1 - \mu} (\dot{p} \partial_y + \dot{q} \partial_y) \dot{u}\bigr) = 0 ,
 \end{gathered}
}
where $p, \tilde{p}, \dot{p} \ge 0$ and $q, \tilde{q}, \dot{q} \in \R$, and $|\dot{q}| \le \dot{p}$, are Fourier multipliers with symbol
\formula{
 -\sym(\xi) & = -(A + i B \sign \xi) |\xi|^\mu ,
}
where $A \ge 0$, $B \in \R$ and $|\arg(A + i B)| \le \min\{\mu, 2 - \mu\}$, and can be represented as
\formula{
 \dn f(x) & = \left\{\begin{array}{lr}
  \displaystyle \frac{1}{|\Gamma(-\mu)|} \int_{-\infty}^\infty \frac{f(x + z) - f(x)}{|z|^{1 + \mu}} \, (C_+ \ind_{(0, \infty)}(z) + C_- \ind_{(-\infty, 0)}(z)) dz \hspace*{-6em} & \text{\raisebox{-2em}{if $\mu \in (0, 1)$,}} \\[2em]
  \displaystyle \frac{A}{\pi} \int_{-\infty}^\infty \frac{f(x + z) - f(x) - z f'(x) \ind_{(-1, 1)}(z)}{|z|^2} \, dz - B f'(x) \hspace*{-6em} & \text{\raisebox{-2em}{if $\mu = 1$,}} \\[2em]
  \displaystyle \frac{1}{\Gamma(-\mu)} \int_{-\infty}^\infty \frac{f(x + z) - f(x) - z f'(x)}{|z|^{1 + \mu}} \, (C_+ \ind_{(0, \infty)}(z) + C_- \ind_{(-\infty, 0)}(z)) dz \hspace*{-6em} & \text{\raisebox{-2em}{if $\mu \in (1, 2)$,}}
 \end{array}\right.
}
where $C_+, C_- \ge 0$. More precisely, the elliptic equations are all equivalent if
\formula{
 \tilde{p} & = p , & \qquad \tilde{q} & = \frac{1 - \mu}{\mu} \, q , \\
 \dot{p} & = \mu^{\mu - 1} (p^2 + q^2)^{\mu/2} , & \qquad \dot{q} & = \mu^{\mu - 1} q (p^2 + q^2)^{-1/2 + \mu/2} = \frac{q}{\sqrt{p^2 + q^2}} \, \dot{p} ;
}
note that when $\mu = 1$, then $\tilde{q}$ is always $0$, see Section~\ref{sec:const} for further discussion. The corresponding coefficients $A$ and $B$ are given by
\formula{
 A + i B & = \begin{cases}
  \displaystyle \frac{(-\Gamma(-\mu)) \Gamma(\mu + \tfrac{(1 - \mu) (p - i q)}{2 p})}{\Gamma(\mu) \Gamma(\tfrac{(1 - \mu) (p - i q)}{2 p})} \, (2 \mu p)^\mu & \text{if $\mu \ne 1$ and $p > 0$,} \\[1.5em]
  \displaystyle \frac{-\Gamma(-\mu)}{\Gamma(\mu)} \, e^{-(i \pi \mu/2) \sign (q (1 - \mu))} |q \mu (1 - \mu)|^\mu & \text{if $\mu \ne 1$ and $p = 0$,} \\[1em]
  \displaystyle p - i q & \text{if $\mu = 1$.}
 \end{cases}
}
Finally, when $\mu \ne 1$, the relation between $(A, B)$ and $(C_+, C_-)$ is determined by
\formula{
 \begin{gathered}
  A + i B = (C_+ e^{-i \mu \pi/2} + C_- e^{i \mu \pi/2}) \sign(1 - \mu) , \\
  \begin{aligned}
   C_+ & = \biggl\lvert \frac{A}{2 \cos \tfrac{\mu \pi}{2}} - \frac{B}{2 \sin \tfrac{\mu \pi}{2}} \biggr\rvert , &  \qquad C_- & = \biggl\lvert \frac{A}{2 \cos \tfrac{\mu \pi}{2}} + \frac{B}{2 \sin \tfrac{\mu \pi}{2}} \biggr\rvert .
  \end{aligned}
 \end{gathered}
}

%
%                            ---------- o ----------
%

\appendix

\section{Proof of the direct part of the representation theorem}
\label{app:proof}

In Section~\ref{sec:ode} we discussed the properties of solutions of the second-order ordinary differential equation
\formula[eqa:ode]{
 \ph'' & = \xi^2 a \ph - 2 i \xi b \ph'
}
(see~\eqref{eq:ode}). Here $\ph$ is assumed to be a continuous function on $[0, R)$ such that the second distributional derivative $\ph''$ corresponds to a measure. In this case necessarily $\ph$ is absolutely continuous, and the distributional derivative $\ph'$ corresponds to a function of bounded variation, equal almost everywhere to the pointwise derivative of $\ph$. Throughout this section, we denote by $\ph'(t)$ the left-continuous version of the point-wise derivative of $\ph$. Note that with this convention, if $\ph$ is a solution of~\eqref{eqa:ode}, then $\ph'(0^+) - \ph'(0) = \xi^2 a(\{0\}) \ph(0)$.

Unlike in Section~\ref{sec:ode}, here we omit the arguments of functions and measures whenever this causes no confusion. For example, we write equations as in~\eqref{eqa:ode} rather than as in~\eqref{eq:ode}.

For a given $\xi > 0$, our goal is to construct a solution $\ph$ of~\eqref{eq:ode} such that $\ph(0) = 1$ and either $\ph$ is a bounded function on $[0, R)$ (if $R = \infty$) or $\ph(R^-) = 0$ (if $R < \infty$). We also need to prove various properties of this solution; most notably, that the mapping $\xi \mapsto -\ph'(0)$ extends to a Rogers function of $\xi$.

We divide the argument into a number of lemmas. The first one is a completely standard application of Picard's iteration. For the convenience of the reader, we provide full details.

\begin{lemma}
\label{lem:picard}
The space of solutions of~\eqref{eqa:ode} is spanned by two linearly independent solutions $\ph_D$ and $\ph_N$, satisfying the initial conditions $\ph_D(0) = \ph_N'(0) = 0$, $\ph_D'(0) = \ph_N(0) = 1$. Furthermore, for every $t \in [0, R)$ the values $\ph_D(t)$, $\ph_D'(t)$, $\ph_N(t)$ and $\ph_N'(t)$ are entire functions of $\xi$.
\end{lemma}

\begin{proof}
Clearly, $\ph$ is a solution of~\eqref{eqa:ode} with initial conditions $\ph(0) = \alpha$, $\ph'(0) = \beta$ if and only if for $t \in [0, R)$ we have
\formula[eqa:ode:int]{
 \ph'(t) & = \beta + \xi^2 \int_{[0, t)} \ph(s) a(ds) - 2 i \xi \int_0^t b(s) \ph'(s) ds , \\
 \ph(t) & = \alpha + \int_0^t \ph'(s) ds .
}
Existence of the solution of~\eqref{eqa:ode:int} on $[0, R)$ follows by Banach's fixed point theorem. In order to define an appropriate Banach space, we choose $C \ge |\xi|$ and we introduce an auxiliary function $M$, defined by
\formula{
 M(t) & = \exp\biggl(2 t + 4 C^2 a([0, t)) + 8 C \int_0^t |b(s)| ds \biggr) .
}
It is easy to see that
\formula{
 M(t) \ge M(0) & = 1 , & \int_0^t M(s) a(ds) & \le \frac{M(t)}{8 C} \, , \\
 M(0) + \int_0^t M(s) ds & \le \frac{M(t)}{2} \, , & \int_0^t M(s) |b(s)| ds & \le \frac{M(t)}{4 C^2} \, .
}
We now consider the Banach space $X$ of absolutely continuous functions $\ph$ such that the second distributional derivative $\ph''$ corresponds to a measure, and the norm in $X$, defined by
\formula{
 \|\ph\|_X & = |\ph(0)| + \sup \biggl\{\frac{|\ph'(t)|}{M(t)} : t \in [0, R]\biggr\} ,
}
is finite. Here, as usual, $\ph'$ corresponds to the left-continuous version of the derivative of $\ph$. Observe that if $\ph \in X$, then
\formula{
 |\ph(t)| & \le |\ph(0)| + \int_0^t |\ph'(t)| dt \\
 & \le |\ph(0)| + \biggl(\sup \biggl\{\frac{|\ph'(s)|}{M(s)} : s \in [0, t]\biggr\}\biggr) \int_0^t M(s) ds \le \|\ph\|_X M(t)
}
for $t \in [0, R)$. In other words, both $|\ph| / M$ and $|\ph'| / M$ are bounded by $\|\ph\|_X$ on $[0, R)$. Finally, we introduce an integral operator $\ii$ defined by
\formula{
 (\ii \ph)'(t) & = \beta + \xi^2 \int_{[0, t)} \ph(s) a(ds) - 2 i \xi \int_0^t b(s) \ph'(s) ds , \\
 \ii \ph(t) & = \alpha + \int_0^t (\ii \ph)'(s) ds .
}
First of all, $\ii$ is a well-defined operator on $X$: if $\ph \in X$, then $|\ph| / M$ and $|\ph'| / M$ are bounded by $\|\ph\|_X$, and hence
\formula{
 |(\ii \ph)'(t)| & \le |\beta| + C^2 \|\ph\|_X \int_{[0, t)} M(s) a(ds) + 2 C \|\ph\|_X \int_0^t |b(s)| M(s) ds \\
 & \le |\beta| M(t) + \frac{M(t) \|\ph\|_X}{4} + \frac{M(t) \|\ph\|_X}{4} \, ,
}
and consequently $\|\ii \ph\|_X \le |\alpha| + |\beta| + \tfrac{1}{2} \|\ph\|_X$. In particular, indeed $\ii \ph$ belongs to $X$. In a similar way, if $\ph_1, \ph_2 \in X$, then
\formula{
 \frac{|(\ii \ph_1)'(t) - (\ii \ph_2)'(t)|}{\|\ph_1 - \ph_2\|_X} & \le C^2 \int_{[0, t)} M(s) a(ds) + 2 C \int_0^t |b(s)| M(s) ds \le \frac{M(t)}{4} + \frac{M(t)}{4} \, ,
}
and therefore
\formula{
 \|\ii \ph_1 - \ii \ph_2\|_X & \le \frac{\|\ph_1 - \ph_2\|_X}{2} \, .
}
It follows that $\ii$ is a contraction on $X$, and thus, by Banach's fixed point theorem, $\ii$ has a unique fixed point $\ph$ in $X$. By definition, $\ph(0) = \ii \ph(0) = \alpha$ and $\ph'(0) = (\ii \ph)'(0) = \beta$, and since $\ph = \ii \ph$, we conclude that $\ph$ is a solution of~\eqref{eqa:ode:int} with the desired initial conditions.

In addition, $\ph$ is the limit in $X$ of the iterates $\ph_n = \ii^n \ph_0$ of $\ii$ applied to $\ph_0(t) = 0$. Observe that
\formula{
 \|\ph - \ph_n\|_X & \le \sum_{j = n}^\infty \|\ph_{j + 1} - \ph_j\|_X \le \|\ph_1 - \ph_0\|_X \sum_{j = n}^\infty 2^{-j} \\
 & = 2^{1 - n} \|\ph_1 - \ph_0\|_X = 2^{-n} \|\ii \ph_1\|_X = 2^{1 - n} (|\alpha| + |\beta|) .
}
Therefore, $\|\ph\|_X$ is uniformly bounded with respect to $\xi$ such that $|\xi| \le C$, and the convergence of $\ph_n$ to $\ph$ in $X$ is uniform in this region. It follows that for every $r \in [0, R)$, $\ph_n(t)$ and $\ph_n'(t)$ are uniformly bounded with respect to $t \in [0, r)$ and $\xi$ such that $|\xi| \le C$, and in this region $\ph_n(t)$ and $\ph_n'(t)$ converge uniformly to $\ph(t)$ and $\ph'(t)$. By Morera's theorem and induction, for every $t \in [0, r)$, $\ph_n(t)$ and $\ph_n'(t)$ are holomorphic functions of $\xi$ in the region $|\xi| < C$, and by Morera's theorem and the dominated convergence theorem, $\ph(t)$ and $\ph'(t)$ have a similar property. Since $C > 0$ and $r \in [0, R)$ are arbitrary, we conclude that $\ph(t)$ and $\ph'(t)$ are entire functions of $\xi$ for every $t \in [0, R)$.

By setting $\alpha = 1$ and $\beta = 0$, we obtain existence of $\ph_N$. Similarly, $\alpha = 0$ and $\beta = 1$ lead to existence of $\ph_D$. Clearly, these functions are linearly independent, and their linear combinations are solutions to~\eqref{eqa:ode}. Furthermore, for every $t \in [0, R)$, $\ph_D(t)$, $\ph_D'(t)$, $\ph_N(t)$ and $\ph_N'(t)$ are entire functions of $\xi$.

Banach's fixed point theorem asserts that $\ph_D$ and $\ph_N$ are unique in $X$. To prove uniqueness of $\ph_D$ and $\ph_N$ in the general class of admissible functions $\ph$, one observes that if $\ph$ is a solution of~\eqref{eqa:ode}, then $\ph'$ is a function with bounded variation, so that $|\ph'| / M$ is bounded on every interval $[0, r)$, where $r \in [0, R)$. Repeating the above proof with $X$ replaced by the Banach space $X_r$ defined in a similar way, but with $R$ replaced by $r$, one obtains uniqueness of solutions on every interval $[0, r)$, with $r \in [0, R)$. Of course this implies that $\ph_D$ and $\ph_N$ are unique solutions on $[0, R)$, and every solution is a linear combination of $\ph_D$ and $\ph_N$.
\end{proof}

The next lemma is a key technical result. Recall that we write $\tilde{a}(dt) = a(dt) - (b(t))^2 dt$ and $B(t) = \int_0^t b(s) ds$.

\begin{lemma}
\label{lem:monotone}
Suppose that $\re \xi > 0$ and $\ph$ is a solution of~\eqref{eqa:ode}. Then
\formula{
 & e^{-2 B \im \xi} \re(\overline{\xi \ph} \ph')
}
is a non-decreasing function on $[0, R)$.
\end{lemma}

\begin{proof}
For $t \in [0, R)$, denote
\formula{
 f(t) & = \overline{\xi \ph(t)} \ph'(t) .
}
Since $\ph$ is a solution of~\eqref{eqa:ode}, the distributional derivative $f'$ corresponds to a measure, which satisfies
\formula{
 f' & = \overline{\xi} (|\ph'|^2 + \overline{\ph} \ph'') \\
 & = \overline{\xi} (|\ph'|^2 + \xi^2 a |\ph|^2 - 2 i \xi b \overline{\ph} \ph') \\
 & = \overline{\xi} |\ph'|^2 + \xi |\xi|^2 a |\ph|^2 - 2 i \xi b f .
}
After elementary manipulations, we find that
\formula{
 \re f' & = (|\xi|^2 a |\ph|^2 + |\ph'|^2 ) \re \xi - 2 \re(i \xi b f) \\
 & = (|\xi|^2 \tilde{a} |\ph|^2 + |\ph' + i \xi b \ph|^2 + 2 b \re(i \overline{\xi \ph} \ph')) \re \xi - 2 b \re(i \xi f) \\
 & = (|\xi|^2 \tilde{a} |\ph|^2 + |\ph' + i \xi b \ph|^2) \re \xi + 2 b \re(i f) \re \xi - 2 b \re(i \xi f) \\
 & = (|\xi|^2 \tilde{a} |\ph|^2 + |\ph' + i \xi b \ph|^2) \re \xi + 2 b \re f \im \xi .
}
Since $(e^{-2 B \im \xi} \re f)' = e^{-2 B \im \xi} (-2 b \re f \im \xi + \re f')$, we find that
\formula[eqa:monotone]{
 (e^{-2 B \im \xi} \re f)' & = e^{-2 B \im \xi} (|\xi|^2 \tilde{a} |\ph|^2 + |\ph' + i \xi b \ph|^2) \re \xi \ge 0 ,
}
that is, $e^{-2 B \im \xi} \re f$ is a non-decreasing function, as desired.
\end{proof}

It is convenient to re-write the assertion of Lemma~\ref{lem:monotone} in terms of the function
\formula{
 \tilde{\ph}(t) & = e^{i \xi B(t)} \ph(t) ,
}
introduced in Section~\ref{sec:ode}, and used frequently below. Since $\tilde{\ph}' = e^{i \xi B} (\ph' + i \xi b \ph)$, we have
\formula{
 |\tilde{\ph}|^2 & = e^{-2 B \im \xi} |\ph|^2 , & |\tilde{\ph}'|^2 & = e^{-2 B \im \xi} |\ph' + i \xi b \ph|^2 , & \re(\overline{\xi \tilde{\ph}} \tilde{\ph}') & = e^{-2 B \im \xi} \re(\overline{\xi \ph} \ph') .
}
Therefore, formula~\eqref{eqa:monotone} reads
\formula[eqa:monotone:alt]{
 (\re(\overline{\xi \tilde{\ph}} \tilde{\ph}'))' & = (|\xi|^2 \tilde{a} |\tilde{\ph}|^2 + |\tilde{\ph}'|^2) \re \xi \ge 0 .
}
Note that although the left-hand side is always a measure, the distributional derivative of $\overline{\xi \tilde{\ph}} \tilde{\ph}'$ (rather than the real part of this function) need not correspond to a measure.

When $\xi > 0$, formulae~\eqref{eqa:monotone} and~\eqref{eqa:monotone:alt} simplify as described in the next result.

\begin{lemma}
\label{lem:convex}
If $\xi > 0$, $\ph$ is a solution of~\eqref{eqa:ode} and $\tilde{\ph}(t) = e^{i \xi B(t)} \ph(t)$ (with $B(t) = \int_0^t b(s) ds$), then $|\ph|^2 = |\tilde{\ph}|^2$ is a convex function on $[0, R)$, and
\formula[eqa:monotone:real]{
 \tfrac{1}{2} (|\ph|^2)' & = \re(\overline{\ph} \ph')' = \xi^2 \tilde{a} |\ph|^2 + |\ph' + i \xi b \ph|^2 = \xi^2 \tilde{a} |\tilde{\ph}|^2 + |\tilde{\ph}'|^2 \ge 0 .
}
Furthermore, if $|\ph|^2$ is non-decreasing, then also $|\ph'|^2$ is non-decreasing, while if $|\ph|^2$ is non-increasing, then also $|\ph'|^2$ is non-increasing.
\end{lemma}

\begin{proof}
The first assertion follows directly from~\eqref{eqa:monotone} and~\eqref{eqa:monotone:alt}. Furthermore, by a direct calculation, 
\formula{
 (|\ph'|^2)' & = 2 \re (\overline{\ph}'' \ph') = 2 \re (\xi^2 a \overline{\ph} \ph' - 2 i \xi b |\ph'|^2) = 2 \xi^2 a \re(\overline{\ph} \ph')
}
has the same sign as $(|\ph|^2)'$.
\end{proof}

\begin{lemma}
\label{lem:unbounded}
If $\xi > 0$ and $\ph_D$, $\ph_N$ are defined as in Lemma~\ref{lem:picard}, then $|\ph_D|^2$ is convex and increasing, while $|\ph_N|^2$ is convex and non-decreasing. If $R = \infty$, then $\ph_D$ is unbounded, and $\ph_N$ is either unbounded or constant.
\end{lemma}

\begin{proof}
Convexity of $|\ph_D|^2$ and $|\ph_N|^2$ is granted by Lemma~\ref{lem:convex}. Since $\ph_D'(0^+) = 1 > 0$, we have $|\ph_D(t)|^2 > 0 = |\ph_D(0)|^2$ for $t \in (0, R)$ small enough. This property and convexity imply that $(|\ph_D|^2)' > 0$ on $(0, R)$, and thus $|\ph_D|^2$ is increasing on $[0, R)$. Furthermore, $\ph_N(0) = 1$ and $\ph_N'(0^+) = \xi^2 a(\{0\}) \ph_N(0) \ge 0$, so that $(|\ph_N|^2)'(0^+) \ge 0$. By convexity, $|\ph_N|^2$ is non-decreasing on $[0, R)$. Finally, a non-decreasing convex function in $[0, \infty)$ is either constant or unbounded.
\end{proof}

%Although we impose a Dirichlet condition on $\ph$ at $R$ whenever $R$ is finite, it turns out that in some cases it is sufficient to require boundedness of $\ph$ near $R$ (and, consequently, boundedness of the $\leb^2(\R)$-valued function $y \mapsto u(\cdot, y)$ in Definition~\ref{def:harm}).

%\begin{lemma}
%\label{lem:bounded}
%If $R < \infty$ and $\xi > 0$, then the function $\ph_D$ defined in Lemma~\ref{lem:picard} is unbounded if and only if [[TODO]].
%\end{lemma}
%
%\begin{proof}
%{}[[TODO]]
%\end{proof}

We now come to the main results of this section, which we split into the following two lemmas.

\begin{lemma}
\label{lem:phi}
If $R = \infty$ and $\xi > 0$, then there is a unique bounded solution $\ph$ of~\eqref{eqa:ode} such that $\ph(0) = 1$, and every other solution diverges to infinity at $\infty$. If $R < \infty$ and $\xi > 0$, then there is a unique solution $\ph$ of~\eqref{eqa:ode} such that $\ph(0) = 1$ and $\ph(R^-) = 0$, and every other solution is bounded away from zero in some left neighbourhood of $R$.
\end{lemma}

\begin{proof}
Let $\ph_D$ and $\ph_N$ be the solutions described in Lemma~\ref{lem:picard}. For $r \in (0, R)$ we define
\formula{
 \beta_r & = -\frac{\ph_N(r)}{\ph_D(r)} && \text{and} & \ph_r & = \ph_N + \beta_r \ph_D ,
}
so that $\ph_r$ is a solution of~\eqref{eqa:ode} satisfying $\ph_r(0) = 1$ and $\ph_r(r) = 0$. Note that $\ph_D(r) \ne 0$, so that $\beta_r$ and $\ph_r$ are well-defined. Our goal is to prove that $\ph = \lim_{r \to R^-} \ph_r$ is the desired solution of~\eqref{eq:ode}.

Suppose that $R = \infty$. Since $|\ph_r|^2$ is convex by Lemma~\ref{lem:convex}, we have $|\ph_r(t)| \le \max\{|\ph_r(0)|, |\ph_r(r)|\} = 1$ for $t \in [0, r]$. It follows that if $0 < t < r$, then
\formula{
 |\beta_t - \beta_r| & = \frac{|\ph_t(t) - \ph_r(t)|}{|\ph_D(t)|} \le \frac{2}{|\ph_D(t)|} \, .
}
By Lemma~\ref{lem:unbounded}, $\lim_{t \to \infty} |\ph_D(t)| = \infty$, so that $|\beta_t - \beta_r| \to 0$ as $t, r \to \infty$. It follows that a finite limit $\beta = \lim_{r \to \infty} \beta_r$ exists, and if we let
\formula{
 \ph(t) & = \lim_{r \to \infty} \ph_r(t) = \ph_D(t) + \beta \ph_N(t) ,
}
then $\ph$ is a bounded solution of~\eqref{eqa:ode} satisfying $\ph(0) = 1$. Every other solution of~\eqref{eqa:ode} which takes value $1$ at $0$ is given by $\ph + \gamma \ph_D$ for some $\gamma \in \C$. Since $\ph_D$ diverges to infinity at $\infty$, $\ph$ is the unique bounded solution of~\eqref{eq:ode} satisfying $\ph(0) = 1$ follows, and every other solution diverges to infinity at $\infty$.

If $R < \infty$, the argument is very similar. By convexity, $|\ph_r(t)|^2 \le 1 - t / r \le 1 - t / R$ for $t \in [0, r]$, so that if $0 < t < r < R$, then
\formula{
 |\beta_t - \beta_r| & = \frac{|\ph_t(t) - \ph_r(t)|}{|\ph_D(t)|} \le \frac{2 \sqrt{1 - t / R}}{|\ph_D(t)|} \, .
}
Since $|\ph_D(t)|$ is increasing, again $|\beta_t - \beta_r| \to 0$ as $t, r \to R^-$, and thus a finite limit $\beta = \lim_{r \to R^-} \beta_r$ exists. We let
\formula{
 \ph(t) & = \lim_{r \to R^-} \ph_r(t) = \ph_D(t) + \beta \ph_N(t) .
}
Clearly, $\ph$ is a solution of~\eqref{eqa:ode} satisfying $\ph(0) = 1$, and since $|\ph(t)|^2 \le 1 - t / R$ for $t \in [0, R)$, we also have $\ph(R^-) = 0$. Finally, since $|\ph_D|$ is increasing, $\ph_D$ is bounded away from zero in some left neighbourhood of $R$. Thus, $\ph$ is the unique solution of~\eqref{eqa:ode} such that $\ph(0) = 1$ and $\ph(R^-) = 0$, and every other solution is bounded away from zero near $R^-$.
\end{proof}

\begin{lemma}
\label{lem:phiprop}
The solution $\ph$ of~\eqref{eqa:ode} described in Lemma~\ref{lem:phi} has the following properties:
\begin{enumerate}[label={\textnormal{(\alph*)}}]
\item\label{lem:phiprop:a} for every $\xi > 0$, $|\ph|^2$ is positive, non-increasing and convex on $[0, R)$, and $|\ph'|$ is non-increasing on $[0, R)$;
\item\label{lem:phiprop:b} the function $\xi \mapsto -\ph'(0)$ extends to a Rogers function $\sym$ (that is, to a holomorphic function $\sym$ in the right complex half-plane, satisfying $\re(\sym(\xi) / \xi) \ge 0$ for every $\xi \in \C$ such that $\re \xi > 0$);
\item\label{lem:phiprop:c} if $\xi > 0$ and $\ph_D$ and $\ph_N$ are the solutions of~\eqref{eqa:ode} described in Lemma~\ref{lem:picard}, then
\formula{
 \sym(\xi) & = \lim_{t \to R^-} \frac{\ph_N(t)}{\ph_D(t)} = \biggl(\int_0^R \frac{e^{-2 i \xi B(t)}}{(\ph_N(t))^2} \, dt\biggr)^{-1} ,
}
where $B(t) = \int_0^t b(s) ds$;
\item\label{lem:phiprop:d} if $\re \xi > 0$, $\tilde{a}(dt) = a(dt) - (b(t))^2 dt$, $B(t) = \int_0^t b(s) ds$ and $\tilde{\ph}(t) = e^{i \xi B(t)} \ph(t)$, then
\formula[eq:phiprop:d]{
 |\xi|^2 \int_{[0, R)} |\tilde{\ph}(t)|^2 \tilde{a}(dt) + \int_0^R |\tilde{\ph}'(t)|^2 dt & = \frac{|\xi|^2}{\re \xi} \, \re \frac{\sym(\xi)}{\xi} \, ,
}
and if $\xi > 0$ and $t \in (0, R)$, then additionally
\formula{
 \xi^2 \int_{[t, R)} |\tilde{\ph}(s)|^2 \tilde{a}(ds) + \int_t^R |\tilde{\ph}'(s)|^2 ds & \le \frac{1}{2 t} \, .
}
\end{enumerate}
In item~\ref{lem:phiprop:d}, for a fixed $t \in [0, R)$, $\ph(t)$ and $\ph'(t)$ denote the holomorphic extensions of functions $\xi \mapsto \ph(t)$ and $\xi \mapsto \ph'(t)$, initially defined for $\xi > 0$.
\end{lemma}

\begin{proof}
The function $|\ph|^2$ is convex by Lemma~\ref{lem:convex}. A convex function on $[0, \infty)$ is either non-increasing or unbounded. If $R = \infty$, then $|\ph|^2$ is bounded, and hence it is non-increasing. When $R < \infty$, then $|\ph|^2$ is non-negative, convex, and it converges to zero at $R^-$. Again, this implies that $|\ph|^2$ is non-increasing. By Lemma~\ref{lem:convex}, also $|\ph'|^2$ is non-increasing.

In order to complete the proof of part~\ref{lem:phiprop:a}, observe that if $\ph(r) = 0$ for some $r \in [0, R)$, then monotonicity of $|\ph|^2$ implies that $\ph(t) = 0$ for all $t \in [r, R)$, and hence, by uniqueness of solutions, $\ph(t) = 0$ for all $t \in [0, R)$, a contradiction. Thus indeed $\ph \ne 0$ on $[0, R)$.

For the proof of part~\ref{lem:phiprop:b}, we use the notation $\ph_D$ and $\ph_N$ introduced in Lemma~\ref{lem:picard}, and we define $\ph_r = \ph_N + \beta_r \ph_D$, where $\beta_r = -\ph_N(r) / \ph_D(r)$, as in the proof of Lemma~\ref{lem:phi}, but for a general $\xi \in \C$ such that $\re \xi > 0$. By Lemma~\ref{lem:monotone}, $\re(\overline{\xi \ph_r(0)} \ph_r'(0)) \le 0$. Since $\ph_r(0) = 1$ and $\ph_r'(0) = \beta_r$, we have $\re(\beta_r / \xi) = |\xi|^{-2} \re(\overline{\xi} \ph_r'(0)) \le 0$. It follows that the mapping $\xi \mapsto -\beta_r$ is a Rogers function. It remains to note $\lim_{r \to R^-} \ph_r'(0) = \lim_{r \to R^-} \beta_r = \beta = \ph'(0)$ whenever $\xi > 0$, and a point-wise limit of Rogers functions on $(0, \infty)$ necessarily extends to a Rogers function (see Remark~3.16 in~\cite{k19}).

We proceed to the proof of part~\ref{lem:phiprop:c}. If $R = \infty$ and $\xi > 0$, then $\ph$ is bounded and $\lim_{t \to \infty} |\ph_D(t)| = \infty$ by Lemma~\ref{lem:unbounded}. Therefore,
\formula{
 \lim_{t \to \infty} \frac{\ph_N(t)}{\ph_D(t)} + \ph'(0) & = \lim_{t \to \infty} \frac{\ph(t)}{\ph_D(t)} = 0 .
}
Similarly, if $R < \infty$, then
\formula{
 \lim_{t \to R^-} \frac{\ph_N(t)}{\ph_D(t)} + \ph'(0) & = \lim_{t \to R^-} \frac{\ph(t)}{\ph_D(t)} = 0 ,
}
because $\ph(R^-) = 0$ and $|\ph_D(t)|$ is increasing by Lemma~\ref{lem:unbounded}. Furthermore, the Wrońskian $W = \ph_D' \ph_N - \ph_D \ph_N'$ satisfies $W(0) = 1$ and $W' = -2 i \xi b W$, so that $W(t) = e^{-2 i \xi B}$. Thus,
\formula{
 \int_0^R \frac{e^{-2 i \xi B(t)}}{(\ph_N(t))^2} \, dt & = \int_0^R \biggl(\frac{\ph_D}{\ph_N}\biggr)'\!(t) dt = \lim_{t \to R^-} \frac{\ph_D(t)}{\ph_N(t)} - \frac{\ph_D(0)}{\ph_N(0)} = \frac{1}{-\ph'(0)} = \frac{1}{\sym(\xi)} \, ,
}
as desired (here we understand that $1 / 0 = \infty$ if $\ph'(0) = 0$).

In order to prove part~\ref{lem:phiprop:d}, we fix $\xi$ such that $\re \xi > 0$, and we again use the notation $\ph_r$ already introduced above. We also write $\tilde{\ph} = e^{i \xi B} \ph$, and similarly we let $\tilde{\ph}_r = e^{i \xi B} \ph_r$. By Lemma~\ref{lem:monotone}, $\re(\overline{\xi \tilde{\ph}_r} \tilde{\ph}_r')$ is non-decreasing (see~\eqref{eqa:monotone:alt}), and clearly $\tilde{\ph}_r(r) = 0$, so that $\re(\overline{\xi \tilde{\ph}_r} \tilde{\ph}_r') \le 0$ on $[0, r]$. Passing to the limit as $r \to R^-$, we find that $\re(\overline{\xi \tilde{\ph}} \tilde{\ph}')$ is non-decreasing and non-positive on $[0, R)$. Furthermore, by~\eqref{eqa:monotone:alt},
\formula[eq:phiprop:d:aux]{
 \lim_{t \to R^-} \re(\overline{\xi \tilde{\ph}(t)} \tilde{\ph}'(t)) - \re(\overline{\xi \tilde{\ph}(0)} \tilde{\ph}'(0)) & = \int_{[0, R)} (\re(\overline{\xi \tilde{\ph}} \tilde{\ph}'))' \\
 & \hspace*{-7em} = |\xi|^2 \re \xi \int_{[0, R)} \tilde{a} |\tilde{\ph}|^2 + \re \xi \int_0^R |\tilde{\ph}'|^2
}
(where for simplicity we omit the argument in the integrands). Since $\re(\overline{\xi \tilde{\ph}} \tilde{\ph}') = \re(\overline{\xi \ph} (\ph' + i \xi b \ph)) = \re(\overline{\xi \ph} \ph')$, we have
\formula{
 \re(\overline{\xi \tilde{\ph}(0)} \tilde{\ph}'(0)) & = \re(\overline{\xi \ph(0)} \ph'(0)) = -\re(\overline{\xi} \sym(\xi)) .
}
We claim that the limit in the left-hand side of~\eqref{eq:phiprop:d:aux} is equal to zero. Together with the above equality and~\eqref{eq:phiprop:d:aux}, this will lead to~\eqref{eq:phiprop:d}.

In order to prove our claim, observe that the non-decreasing, non-positive function $\re(\overline{\xi \tilde{\ph}} \tilde{\ph}')$ necessarily has a finite limit at $R^-$. Suppose, contrary to our claim, that this limit is non-zero. Then there is $C_1 > 0$ such that $|\overline{\tilde{\ph}} \tilde{\ph}'| \ge C_1$ in some left neighbourhood of $R$. We now consider two cases. If $R = \infty$, then, by Schwarz inequality, $|\tilde{\ph}(t)| \le 1 + \sqrt{C_2 t}$, where $C_2$ is the integral of $|\tilde{\ph}'|^2$. Therefore, $|\tilde{\ph}'|^2 \ge C_1^2 (1 + \sqrt{C_2 t})^{-2}$ in some neighbourhood of $\infty$, and hence $|\tilde{\ph}'|^2$ is not integrable. This is a contradiction: by~\eqref{eq:phiprop:d}, $|\tilde{\ph}'|^2$ is integrable. It follows that $\re(\overline{\xi \tilde{\ph}} \tilde{\ph}')$ indeed converges to zero at $\infty$, as desired.

If $R < \infty$, the argument slightly more involved. Recall that $\tilde{\ph}$ is the limit of $\tilde{\ph}_r$ as $r \to R^-$, and by~\eqref{eqa:monotone:alt}, for every $r \in (0, R)$ we have
\formula{
 -\re(\overline{\xi \tilde{\ph}_r(0)} \tilde{\ph}_r'(0)) & = \int_{[0, r)} (\re(\overline{\xi \tilde{\ph}_r} \tilde{\ph}_r'))' = |\xi|^2 \re \xi \int_{[0, r)} \tilde{a} |\tilde{\ph}_r|^2 + \re \xi \int_0^r |\tilde{\ph}_r'|^2 ,
}
just as in~\eqref{eq:phiprop:d:aux}. Since $\overline{\xi \tilde{\ph}_r(0)} \tilde{\ph}_r'(0)$ converges to a finite limit $\overline{\xi \tilde{\ph}(0)} \tilde{\ph}'(0)$ as $r \to R^-$, it follows that $\int_0^r |\tilde{\ph}_r'|^2 \le C_3$ for every $r \in (R/2, R)$ and some $C_3 > 0$. By Schwarz inequality, we have $|\tilde{\ph}_r(t)| \le \sqrt{C_3 (r - t)} \le \sqrt{C_3 (R - t)}$ for $t \in [0, r]$. Passing to the limit as $r \to R^-$, we find that $|\tilde{\ph}(t)| \le \sqrt{C_3 (R - t)}$ for $t \in [0, R]$, and therefore $|\tilde{\ph}'|^2 \ge C_1^2 C_3^{-1} (R - t)^{-1}$. This again contradicts integrability of $|\tilde{\ph}'|^2$ asserted by~\eqref{eq:phiprop:d}, and our claim follows.

It remains to prove the other part of item~\ref{lem:phiprop:d} of the lemma. Observe that if $\xi > 0$, then $|\tilde{\ph}|^2 = |\ph|^2$ is convex by Lemma~\ref{lem:convex}, and hence
\formula{
 (|\tilde{\ph}|^2)'(t) & \ge \frac{|\tilde{\ph}(t)|^2 - |\tilde{\ph}(0)|^2}{t} \ge -\frac{1}{t} \, .
}
Thus, as in formula~\eqref{eq:phiprop:d:aux}, we have
\formula{
 |\xi|^2 \int_{[t, R)} \tilde{a} |\tilde{\ph}|^2 + \int_t^R |\tilde{\ph}'|^2 & = \int_{[t, R)} (\re(\overline{\tilde{\ph}} \tilde{\ph}'))' = -\re(\overline{\tilde{\ph}(t)} \tilde{\ph}'(t)) = -\frac{(|\tilde{\ph}|^2)'(t)}{2} \le \frac{1}{2 t} \, ,
}
as desired.
\end{proof}

Part~\ref{thm:ode:a} of Theorem~\ref{thm:ode} is now a consequence of Lemma~\ref{lem:phi}, while parts~\ref{thm:ode:b} and~\ref{thm:ode:c} follow from Lemma~\ref{lem:phiprop}.

%
%                            ---------- o ----------
%

\medskip

\subsection*{Acknowledgments}

I thank Tadeusz Kulczycki and Jacek Mucha for valuable discussions on the subject of this article. I thank Tomasz Grzywny and Pablo Raúl Stinga for their comments to the preliminary version of this article.

%
%                            ---------- o ----------
%

%
%                            ---------- o ----------
%

\end{document}